\documentclass[12pt]{article}

\usepackage{amsmath, amsthm, amssymb, amsfonts, enumerate, dsfont}
\usepackage[applemac]{inputenc} 
 \usepackage[french] {babel}

\def\dis{\displaystyle}
\def \R{I\!\!R}
\def \E{I\!\!E}
\def \P{I\!\!P}
\def \N{I\!\!N}

\newtheorem{thm}{Theorem}[section]
\newtheorem{cor}[thm]{Corollary}
\newtheorem{lem}[thm]{Lemma}
\newtheorem{pro}[thm]{Proposition}
\newtheorem{defi}[thm]{Definition}
\newtheorem{rem}[thm]{Remark}

\newtheorem{exa}[thm]{Example}

\newtheorem{ass}[thm]{Assumption}

\def \conv{ {\rm conv} }

\def \cl { {\rm cl} }
\def \conv{ {\rm conv} }
\def \cav{ {\rm cav} }

\def \Min{ {\rm Min} }

\def \Val{ {\rm Val }}

\def \lim{ {\rm lim} }

\def \1{\mathbf{1}}
 \def \C{\mathds{C}}

\newcommand{\Tau}{\text{\Large $\tau$}}

\def\abstract{\begin{center} \small\bf Abstract\end{center}\small}

 \title{A distance for probability spaces, and  long-term  \\ values
  in Markov Decision Processes and Repeated Games}

 \author{J\'er\^ome Renault\thanks{%
TSE (GREMAQ, Universit\' e Toulouse 1 Capitole), 21 all\' ee de Brienne, 31000
Toulouse, France. E-mail: \textsf{jerome.renault@tse-fr.eu}.} , Xavier Venel%
\thanks{%
TSE (GREMAQ, Universit\' e Toulouse 1 Capitole), 21 all\' ee de Brienne, 31000
Toulouse, France.  .}}
\date{\today}

\begin{document}

\maketitle

\begin{abstract}
Given a finite set $K$, we denote by $X=\Delta(K)$ the set of probabilities on $K$ and by $Z=\Delta_f(X)$ the set of Borel probabilities on $X$  with finite support. Studying    a Markov Decision Process with partial information on $K$ naturally 	 leads to a Markov Decision Process with full information on $X$. We introduce a new metric $d_*$ on $Z$ such that  the  transitions become $1$-Lipschitz from $(X, \|.\|_1)$  to $(Z,d_*)$. In the first part of the article, we define and  prove  several properties of  the metric $d_*$. Especially, $d_*$ satisfies a Kantorovich-Rubinstein type  duality formula    and can be characterized by using disintegrations. In the second part, we characterize the limit values in several classes of ``compact non expansive" Markov Decision Processes.  In particular we use the  metric $d_*$  to characterize the limit value in Partial Observation MDP with finitely many states and in Repeated Games with an informed controller with finite sets of states and actions. Moreover in each case we can prove the existence of a  generalized notion of uniform value where we consider not only the Cesàro mean when the number of stages  is large  enough but any evaluation function $\theta \in \Delta(\N^*)$ when the impatience $I(\theta)=\sum_{t\geq 1} |\theta_{t+1}-\theta_t|$ is small enough. \\

    \end{abstract}

\section{Introduction}

The classic model of Markov Decision Processes with finitely many states, particular class of the model of Stochastic Games introduced by Shapley (1953),  was  explicitly introduced by Bellman (1957) in the 1950s and has been extensively studied since then. When the set of  actions is also  finite, Blackwell (1962) proved the existence of a strategy which is optimal for all discount factors close to $0$.%
This model was generalized later to MDPs with Partial Observations (POMDP), (for references see Araposthatis {\it  et al.}  (1993)). The decision maker observes neither the state nor his payoff. Instead at each stage, he receives a signal which depends on the previous state and his previous action. In order to solve this problem a classic approach is to go back to the classic model of MDPs by introducing an auxiliary problem with full observation and Borel state space: the space of belief on the state as shown in Astrom, K.J. (1965), Sawaragi and Yoshikawa (1970) and Rhenius (1974). For optimality criteria like the Cesàro mean and the Abel mean, these two problems are equivalent and the question of the existence of the limit value is the same. Then given some sufficient conditions of ergodicity, one can search for a solution of the Average Cost Optimality Criterion in order to find ``the" value of the MDP,  for example as in Runggaldier and Stettner(1991) or as in Borkar (2000,2007). An introduction to the ACOE in the framework of MDP and the reduction of POMDP can be found in Hern\'{a}ndez-Lerma (1989).  From another point of view, if we know that the limit value exists, the ACOE may be used as a characterization of the value. For finite MDP,  for example, Denardo and Fox (1968) proved  that the limit value is the solution of a linear programming problem deduced from the ACOE. Moreover by standard linear programming results, it is also equal to the solution of a dual problem from which Hordjik and Kallenberg (1979) deduced an optimal strategy. This dual problem focuses on the maximal payoff that the decision maker can guarantee on invariant measures. This approach was extended to different criteria (see Kallenberg 1994 ) and to a convex analytic approach by Borkar ( for references see Borkar 2002) in order to study problems with a countable state space and a compact action space. 

Given an initial POMDP on a finite space $K$, we will follow the usual approach and introduce a $MDP$ on $X=\Delta(K)$ but instead of assuming some ergodicity on the process we will use the structure of $\Delta(K)$ and a new metric on $Z=\Delta_f(\Delta(K))$. We extend and relax the MDP on $\Delta_f(X)$ with a uniformly continuous affine payoff function and non-expansive affine transitions. The structure of $Z$ was already used in Rosenberg, {\it  et al.}  (2002) and in Renault (2011). Under our new metric, we highlight a stronger property since the transitions became $1$-Lipschitz on $Z$ and $Z$ is still precompact. We use this   property  to focus on general evaluations. Given a probability distribution $\theta$ on positive integers, we evaluate a sequence of payoffs $g=(g_t)_{t\geq 1}$ by $\gamma_{\theta}(g)=\sum_t \theta_t g_t$. In a MDP or a POMDP, the $\theta$-value is then defined as the maximum expected payoff that the player can guarantee with this evaluation. Most of the literature   focuses on the $n$-stage game where we consider the Cesàro mean of length $n$, and on the $\lambda$ discounted games, where we consider the Abel mean with parameter $\lambda$. The first type of results focuses on the limit when $n$ converges to $+\infty$ and when $\lambda$ converges to $0$ or the relation between them. When there is no player,  the relation between them is directly linked to a Hardy-Littlewood theorem (see Filar and Sznajder, 1992). One of the limit exists if and only if the other exists and whenever they exist they are equal. Lehrer and Sorin (1992) proved that this result extends to the case where there is one player provided we ask  for uniform convergence. The other approach focuses  on the existence of a good strategy in any long game or for any discount factor close to $0$. We say that the MDP has a uniform value. For MDP with finitely many  states, Blackwell's result (1962) solved both problems. In POMDPs, Rosenberg, {\it  et al.}  (2002) proved the existence of the uniform value when the sets of states, actions and signals are finite, and Renault (2011) removed the finiteness assumption on signals  and actions.

Concerning stochastic games, Mertens and Neyman (1981) proved the existence of the uniform value when the set of states and the set of actions are finite. The model also generalizes to partial information but the existence of possible private information implies a more complex structure on the auxiliary state space. Mertens and Zamir (1985) and Mertens (1987) introduced the universal belief  space  which synthesizes all the information for both players in a general repeated game: their beliefs about the state, their beliefs about the beliefs of the other player, etc... So far, the results always concern  some subclasses of games where we can explicitly write the auxiliary game in a ``small"  tractable set. A lot of work has been done on games with one fully informed player and one player with partial information, introduced by Aumann and Maschler (see reference from 1995). A state is chosen at stage $0$ and remains  fixed for the rest of the game. Renault (2006)  extended the analysis to a general underlying Markov chain on the state space (see also Neyman, 2008). Rosenberg {\it  et al.}  (2004) and Renault (2012a) proved the existence of the uniform value when the informed player can additionally control the evolution of the state variable.\\

The first section is dedicated to the description of  the (pseudo)-distance $d_*$ on $\Delta(X)$ in the general framework when  $X$ is a compact  subset of a normed vector space. We provide different definitions and show that they all define this pseudo-distance. Then we focus on  the case where $X$ is a simplex. We prove that $d_*$ is a real metric and prove a ``Kantorovich-Rubinstein like " duality formula for probabilities with finite support on $X$. We give new definitions and a characterization by the disintegration mapping. The second section focuses on Gambling Houses and standard Markov Decision Processes. We first introduce the definitions of  general limit value and   general uniform value. Then we give sufficient conditions for the existence of the general uniform value and a characterization in several ``compact" cases of Gambling Houses and Markov Decision Processes, including the finite state case. We   study  the limit value as a linear function of the initial probability so there are similarities  with the convex analytic approach,  but we are able to avoid any  assumption on the set of actions. Moreover the MDPs that we are considering may not have $0$-optimal strategies as shown in Renault (2011).  Finally we apply these results to prove  the existence of the general uniform value in finite state POMDPs and repeated games with an informed controller.

\section{A  distance for   belief spaces} \label{sec2}

 \subsection{A pseudo-distance for  probabilities on a compact subset of a normed vector space}
 
 We fix a compact subset $X$ of a real  normed vector space $V$. We denote by $E={\cal C}(X)$ the set of continuous functions from  $X$ to the reals, and by $E_1$ the set of 1-Lipschitz functions in $E$. We denote by  $\Delta(X)$ the set of Borel probability measures on $X$, and for each $x$ in $X$ we write  $\delta_x$ for the Dirac probability measure on $x$.   It is well known that  $\Delta(X)$ is a compact set for the weak-* topology, and this topology can be metrizable by the (Wasserstein) Kantorovich-Rubinstein  distance:
 $$\forall u,v \in \Delta(X), \;\; d_{KR}(u,v)=\sup_{f \in E_1}\;  u(f)-v(f).$$
 
We will introduce  a pseudo-distance on $\Delta(X)$, which  is not greater than $d_{KR}$ and in some cases also metrizes the weak-* topology. We start with several definitions, which will turn out  to be equivalent. Let $u$ and $v$ be in $\Delta(X)$. 
 
\begin{defi} \label{def1}  $$d_1(u,v) =\sup_{f \in D_1}\;  u(f)-v(f),$$
 $${\rm{\it  where}}\;\; D_1=\{f \in E, \forall x,y \in X, \forall a, b \geq 0, \; a f(x)-b f(y)\leq \| a x-b y\|\}.$$ \end{defi}
 
Note that any  linear  functional  in $V'$ with norm 1 induces an element  of $D_1$. $d_1$ is a pseudo-distance on $\Delta(X)$, and $d_1(u,v) =\sup_{f \in D_1}\;  | u(f)-v(f)|$, since if $f$ is in $D_1$, $-f$ is also in $D_1$. We also have $D_1\subset E_1$, so that $d_1(u,v)\leq d_{KR}(u,v)$ and the supremum in the definition of $d_1(u,v)$ is achieved.
  
Given $x$ and $y$ in $X$, there exists a linear functional  $f$ in $V'$ with norm 1 such that $f(y-x)=\|y-x\|$. Then the restriction of $f$ to $X$ is in $D_1$ and $d_1(\delta_x, \delta_y)\geq \|x-y\|$. One can easily  deduce that $d_1(\delta_x, \delta_y)= \|x-y\|$ for $x$ and $y$ in $X$. \\

\begin{exa} \label{exemple1} \rm Consider the particular case where $X=[0,1]$ endowed with the usual norm. Then all $f$ in $D_1$ are linear. As a consequence,   $d_1(u,v)=0$ for $u=1/2 \, \delta_0 +1/2 \,  \delta_1$ and $v=\delta_{1/2}$. We do not have the separation property and $d_1$ is not a distance in this case. 

Let us modify the example.  $X$ now is the set of probability distributions over 2 elements, viewed as $X=\{(x, 1-x), x \in [0,1]\}$.  We use the norm $\|.\|_1$ to measure the distance between $(x,1-x)$ and $(y,1-y)$, so that $V=\R^2$ is endowed with $\|(x_1, x_2)-(y_1,y_2)\|=|x_1-y_1|+|x_2-y_2|$.   Consider    $f$ in $E$ such that  $f((x,1-x))=x(1-x)$ for all $x$. $f$ now belongs to $D_1$, and $d_1(u,v)\geq 1/4>0$ for $u=1/2 \, \delta_0 +1/2 \,  \delta_1$ and $v=\delta_{1/2}$.  One can show that $(\Delta(X),d_1)$ is a compact metric space in this case (see proposition   \ref{pro7} later), and for applications  in this paper $d_1$  will be a particularly useful distance whenever $X$ is a simplex $\Delta(K)$ endowed with $\|x-y\|=\sum_{k\in K} |x^k-y^k|$.
\end{exa}

Furthermore it is known that the Kantorovitch Rubinstein metric on $\Delta(X)$ only depends on the restriction of the norm $\|.\|$ on the set $X$. Especially if for all $x,x'\in X$ such that $x\neq x'$, $\|x-x'\|=2$, then for all $u,v \in \Delta(X)$, $d_{KR}(u,v)=\|u-v\|_1$. This is not the case when considering the metric $d_1$. Two norms on $V$ giving the same metric on $X$ may leads to different pseudo-metrics on $\Delta(X)$.  We consider in the next example different norms on the Euclidean space $\R^K.$

\begin{exa} \label{norme}
\rm We consider $V=\R^K$, $X=\{e_1,..,e_K\}$ the set of canonical vectors of $V$ and a norm such that for all $k\neq k'$, $\|e_k-e_{k'} \|=2$. We know that $d_1$ is smaller than the Kantorovitch-Rubinstein metric, so for all $u \in \Delta(X)$ and $v\in \Delta(X)$, we have $d_1(u,v) \leq \|u-v\|_1$.\\

\rm We first consider the particular case of the norm defined by $\|x-y\|=2^{1-\frac{1}{p}} \|x-y\|_p$ where  $\|x-y\|_p= {\left(\sum_{k=1}^K |x_k-y_k|^p \right) }^{1/p}$ is the usual $L^p$-norm on $\R^K$, 
with  $p$ a fixed positive integer.  Given $u,v\in \Delta(X)$, the function $f$ defined by 
\[\forall k\in K \ f(k)= \begin{cases}
                               \begin{matrix}
                               1 &\text{ if } u(k)\geq v(k) \\
                              -1 & \text{ otherwise,} 
                              \end{matrix}
                              \end{cases}\]
satisfies $u(f)-v(f)=\sum_{k\in K} |u(k)-v(k)|=\|u-v\|_1$. Moreover for all  $a\geq 0$, $b\geq 0$ and $k,k'\in K$ such that $k\neq k'$, we have
\[
a f(k)-bf(k') \leq a+b \leq \frac{2}{2^{1/p}} \left( a^p + b^p \right)^{1/p} =\|ae_k-be_{k'}\|,
\]
and $a f(k)-bf(k) \leq|a-b|\leq |a-b|\frac{2}{2^{1/p}}.$ Therefore $f$ is in $D_1$ and $d_1(u,v)= \|u-v\|_1$, independently\footnote{Similarly,     the same result holds for the case $p=+\infty$, i.e. where $\|x-y\|=2  \|x-y\|_\infty$.} of $p$.  \\

\rm Nevertheless the inequality $d_1(u,v) \leq \|u-v\|_1$ may be strict as in the following example. We consider the case $K=3$ and given a vector $(x_1,x_2,x_3)\in \R^3$, we define the norm $\|(x_1,x_2,x_3)\|= \max(|x_1|+|x_2|,2 |x_3|)$, which satisfies $\|e_1-e_2\|=\|e_2-e_3\|=\|e_3-e_2\|=2$. Let $f$ be a function in $D_1$, then we have among others the following constraints:
\begin{center}
\begin{tabular}{rlrl}
&$\forall a,b\geq 0 $& $a f(e_3)-b f(e_1)$ & $\leq \|(-b,0,a)\|= \max(2a,b)$ \\
and & $\forall a \geq 0  $& $a f(e_2) $& $\leq \|(0,a,0)\|= a$. 
\end{tabular}
\end{center}
\noindent Let $u=(0,1/2,1/2)$, $v=(1,0,0)$ and $f\in D_1$, then
\[
u(f)-v(f)=\frac{1}{2} f(e_2) +\frac{1}{2} f(e_3) - f(e_1) \leq \frac{1}{2}+\max(2/2,1)=\frac{3}{2}.
\]
By symmetry of $D_1,$ we deduce that $d_1(u,v) \leq \frac{3}{2}<\|u-v\|_1$. In fact one can show that $d_1(u,v)=\frac{3}{2}$ by  checking that   the function defined by  $f(e_1)=0$, $f(e_2)=1$ and $f(e_3)=2$  is in $D_1$  and  satisfies $u(f)-v(f)=\frac{3}{2}$.
\end{exa}

 \vspace{0,5cm}
 
 We now give other expressions for the pseudo-distance $d_1$.
 
\begin{defi} \label{def2} 
$$d_2(u,v) =\sup_{(f,g) \in D_2}\;  u(f)+v(g),$$
 $${\rm{\it  where}}\;\;D_2=\{(f,g)\in E \times E, \forall x,y \in X, \forall a, b \geq 0,  \; a f(x)+b g(y)\leq \| a x-b y\|\}.$$ \end{defi}
 
  \vspace{0,5cm}
 
 \begin{defi} \label{def2+} 
$$d_2^+(u,v)=\inf_{\varepsilon>0} d_2^{\varepsilon}(u,v), \; {\rm where }\; d_2^{\varepsilon}(u,v)=\sup_{(f,g) \in D_2^\varepsilon}\;  u(f)+v(g)$$
 $${\rm{\it  and}}\;  \forall  \varepsilon>0,\;\; D_2^\varepsilon=\{(f,g)\in E \times E, \forall x,y \in X, \forall a, b \in [0,1],  \; a f(x)+b g(y)\leq \varepsilon + \| a x-b y\|\}.$$ \end{defi}
 
 \vspace{0,5cm}
 
  \begin{defi} \label{def3} $$d_3(u,v)=\inf_{\gamma\in {\cal M}_3(u,v)} \int_{ X^2\times [0,1]^2} \|\lambda x-\mu y\| d\gamma(x,y,\lambda, \mu),$$
where   ${\cal M}_3(u,v)$ is   the set of  finite positive measures  on $X^2\times [0,1]^2$ satisfying for each $f$   in $E$:
$$\dis  \int_{(x,y,\lambda, \mu)\in X^2\times [0,1]^2} \lambda f(x) d\gamma(x,y,\lambda,\mu)=u(f),\; {\rm \it and}\; \int_{(x,y,\lambda, \mu)\in X^2\times [0,1]^2} \mu f(y) d\gamma(x,y,\lambda,\mu)=v(f).$$
 \end{defi}



\vspace{0,5cm}
 
In the next subsection we will prove the following  result.

 \begin{thm} \label{thm1}
 For all $u$ and $v$ in $\Delta(X)$, $d_1(u,v)=d_2(u,v)=d_2^+(u,v)=d_3(u,v).$
 \end{thm}

\subsection{Proof of theorem \ref{thm1}}

The proof is split into several parts.

\begin{pro} \label{pro1}$d_1 =d_2=d_2^+.$
\end{pro}

It is plain that $d_1 \leq d_2  \leq d_2^+ $, so all we have to prove is $d_2^+\leq d_1$. We start with a lemma.

\begin{lem}\label{lem1} Fix $\varepsilon>0$, and  let  $f$ in $E$ be such that: $\forall x \in X$, $\forall a\in [0,1]$, $a f(x)\leq \varepsilon + a \|x\|.$ Define $\hat{f}$  by:
$$\forall y \in X, \; \hat{f}(y)=\inf_{a\in [0,1], b\in (0,1], x\in X} \frac{1}{b} \left(\varepsilon + \|ax-by\|-a f(x)\right).$$
Then for each $y$ in $X$, $-\|y \| \leq \hat{f}(y) \leq -f(y)+ \varepsilon$. Moreover $\hat{f} \in E_1$,   and: $$\forall x \in X, \forall y \in X, \forall a\in [0,1], \forall b\in [0,1], \;   a\hat{f}(x)-b\hat{f}(y)\leq a\varepsilon + \|by-ax\|.$$

\end{lem}
\noindent{\bf Proof of  lemma \ref{lem1}:}  By assumption on $f$, we have for all $y$ in $X$, $a$ in $[0,1]$, $b$ in $(0,1]$, $x$ in $X$: $ \frac{1}{b} \left(\varepsilon + \|ax-by\|-a f(x)\right)\geq \frac{1}{b} \left( -a\|x\| + \|ax-by \|\right)$ $\geq -\|y\|$. In the definition  of $\hat{f}(y)$, considering $a=b=1$ and $x=y$ yields  $\hat{f}(y)\leq -f(y)+ \varepsilon$. 

Fix $x$ and  $y$ in $X$, $a$ and $b$ in $[0,1]$. We have:
\begin{eqnarray*}
a\hat{f}(x)-b\hat{f}(y)& =& a\inf_{a',b',x'}  \frac{1}{b'} \left(\varepsilon + \|a'x'-b'x\|-a' f(x')\right)\\
 & &  - b \inf_{a'',b'',x''}  \frac{1}{b''} \left(\varepsilon + \|a''x''-b''y\|-a'' f(x'')\right).
 \end{eqnarray*}
If $a=0$, then the inequality $\hat{f}(y)\geq -\|y\|$ leads to $-b\hat{f}(y)\leq b\|y\|$. If $b=0$, choose $a'=0$, $b'=1$ and $x'=x$ to get $a\hat{f}(x)\leq a\varepsilon + \|ax\|$.
  
If $ab>0$, given $\eta>0$, choose $a''$, $b''$, $x''$ $\eta$-optimal in the second infimum. We can define $x'=x''$, and choose $a'\in [0,1]$ and $b'\in (0,1]$ such that $\frac{a'}{b'}=\frac{b}{a} \frac{a''}{b''}$. We obtain:
\begin{eqnarray*}
a\hat{f}(x)-b\hat{f}(y)& \leq& b \eta+ (\frac{a}{b'}-\frac{b}{b''}) \varepsilon + (\| \frac{a''}{b''}b x''-a x\|- \|\frac{a''}{b''}b x''- by\|)\\
 & \leq &  b \eta+ (\frac{a}{b'}-\frac{b}{b''}) \varepsilon + \|ax-by \|.
\end{eqnarray*}

If $a=b>0$, choose $a'=a''$ and $b'=b''$ to obtain: $ \hat{f}(x)- \hat{f}(y)\leq \|x-y\|$ and therefore $\hat{f}$ is 1-Lipschitz. 

Otherwise, we distinguish two cases. If $\frac{a}{b} b''\leq 1$, we define $b'=\frac{a}{b} b''$ and $a'=a''$ and we get $ a\hat{f}(x)- b\hat{f}(y)\leq b\eta+\|ax-by\|.$ If $\frac{a}{b} b''> 1$, we define $b'=1$ and $a'=\frac{a'' b}{b'' a} \in [0,1]$ and obtain $ a\hat{f}(x)- b\hat{f}(y)\leq b\eta+ a \varepsilon +   \|ax-by\|.$ Thus for all $\eta>0$, we have
\begin{eqnarray*}
a\hat{f}(x)-b\hat{f}(y)& \leq & b\eta + a\varepsilon + \|ax-by \|,
\end{eqnarray*}
and therefore $a\hat{f}(x)-b\hat{f}(y)\leq  a\varepsilon + \|ax-by\|$.\hfill $\Box$

\vspace{0,5cm}

\noindent{\bf Proof of proposition \ref{pro1}:} Fix $u$ and $v$ in $\Delta(X)$, and consider $\varepsilon>0$. For each $(f,g)$ in $D_2^\varepsilon$, we have $-f+ \varepsilon \geq \hat{f}\geq g$ and  $(f,\hat{f})$ in $D_2^\varepsilon$.    We also have $(\hat{f},f)\in D_2^\varepsilon$ so iterating the construction, we get  $( \hat{f}, \hat{\hat f})\in D_2^\varepsilon$, and $-\hat{f}+ \varepsilon \geq \hat{\hat{f}}\geq f$. 

Now, $u(f)+v(g)  \leq u(\hat{\hat{f}}) + v(\hat{f})\leq -u(\hat{f}) + \varepsilon + v(\hat{f}).$ Hence we have obtained:
$$d_2^\varepsilon(u,v)\leq \varepsilon + \sup_{f \in C_{\varepsilon(u,v)}} -u(f)+v(f),$$ where $C_{\varepsilon(u,v)}$ is the set of functions $f$ in $E_1$
 satisfying: $$ \forall x \in X, \forall y \in X, \forall a\in [0,1], \forall b\in [0,1], \;   a f(x)-bf(y)\leq a\varepsilon + \|ax-by\| \; {\it and} \; f(y) \geq - \|y \| .$$
For each positive $k$, one can choose $f_k$ in $E_1$ achieving the above supremum for $\varepsilon=1/k$. Taking a   limit point of $(f_k)_k$ yields  a function $f$ in $D_1$ such that: $-u(f)+v(f)\geq d_2^+(u,v)$. The function $f^*=-f$ is in $D_1$ and satisfies $u(f^*)-v(f^*)\geq   d_2^+(u,v)$, and the proof of proposition \ref{pro1} is complete. \hfill $\Box$
 
 \begin{pro} \label{pro2}
 $d_2^+\geq d_3.$
 \end{pro}
 
 \noindent{\bf Proof:} The proof is based on (a corollary of) Hahn-Banach theorem. Define:
 $H={\cal C}(X^2\times [0,1]^2)$ and $$L=\{\varphi \in H, \exists f,g \in {\cal C} (X) \; s.t. \; \forall x, y \in X, \forall \lambda, \mu \in [0,1], \varphi(x,y,\lambda, \mu)= \lambda f(x) + \mu g(y)\}.$$
 \noindent $H$ is endowed with the uniform norm and $L$ is a linear subspace of $H$. Note that the unique constant mapping in $L$ is 0. Fix $u$ and $v$ in $\Delta(X)$, and let $r$ be the linear form on $L$ defined by $r(\varphi)=u(f)+v(g)$, where $\varphi(x,y,\lambda, \mu)= \lambda f(x) + \mu g(y)$ for all $x$, $y$, $\lambda$, $\mu$.
 
 Fix now $\varepsilon>0$, and put: $$U_{\varepsilon}=\{\varphi \in H, \forall x, y \in X, \forall \lambda, \mu \in [0,1], \varphi(x,y,\lambda, \mu) \leq \|\lambda x- \mu y\| + \varepsilon\}.$$
 We have: $$\sup_{\varphi \in L \cap U_{\varepsilon}} r(\varphi)= d_2^\varepsilon (u,v).$$
 $U_{\varepsilon}$ is a convex subset of $H$ which is radial at 0, in the sense that: $\forall \varphi \in H$, $\exists \delta>0$ such that $t\varphi \in U_{\varepsilon}$ as soon as $|t|\leq \delta$.  By a corollary of Hahn-Banach theorem (see theorem 6.2.11 p.202 in Dudley, 2002), $r$ can be extended to a linear form on $H$ such that: $$\sup_{\varphi \in  U_{\varepsilon}} r(\varphi)= d_2^\varepsilon (u,v).$$

 Given $\varphi \in H$, we have $\varepsilon \varphi / \|\varphi\|_{\infty}\in U_{\varepsilon}$, which implies that $r(\varphi)\leq \|\varphi\|_{\infty}d_2^\varepsilon (u,v)/ \varepsilon$, so that $r$ belongs to $H'$. And if $\varphi \geq 0$, we have $t\varphi \in U_{\varepsilon}$ if $t\leq 0$, so that $r(\varphi)\geq {d_2^\varepsilon (u,v)/t}$ for all $t\leq 0$ and  $r(\varphi)\geq 0$. By Riesz Theorem, $r$ can be represented by a positive finite measure $\gamma$ on $X^2\times [0,1]^2$. 
 
 Given $f$ in $E$, one can consider $\varphi_f\in L$  defined by $\varphi_f(x,y,\lambda, \mu)=\lambda f(x).$ $r(\varphi=f)=\gamma(\varphi_f)$ gives: $u(f)= \int_{(x,y,\lambda, \mu)\in X^2\times [0,1]^2} \lambda f(x) d\gamma(x,y,\lambda,\mu)$, and similarly $v(f)=  \int_{(x,y,\lambda, \mu)\in X^2\times [0,1]^2} \mu f(y) d\gamma(x,y,\lambda,\mu)$, and we obtain that $\gamma \in {\cal M}_3(u,v).$
 
Because $\gamma\geq 0$,  $\sup_{\varphi \in  U_{\varepsilon}} r(\varphi)= r(\varphi^*)$ where $\varphi^*(x,y,\lambda,\mu)=\|\lambda x- \mu y\| + \varepsilon.$ We get $d_2^\varepsilon (u,v)= \int_{ X^2\times [0,1]^2} \|\lambda x-\mu y\| d\gamma(x,y,\lambda, \mu) + \varepsilon \gamma(X^2\times [0,1]^2)$, so  $$d_2^\varepsilon (u,v)\geq \int_{ X^2\times [0,1]^2} \|\lambda x-\mu y\| d\gamma(x,y,\lambda, \mu)\geq d_3(u,v).$$
\noindent   \; \hfill $\Box$ 

\begin{lem} $d_3\geq d_2.$
\end{lem}
 \noindent {\bf Proof:} Fix $(f,g)\in D_2$ and $\gamma \in  {\cal M}_3(u,v).$
 \begin{eqnarray*}
u(f)+v(g) &=& \int_{X^2\times [0,1]^2} \lambda f(x) d\gamma(x,y,\lambda, \mu) + \int_{X^2\times [0,1]^2} \mu g(y) d\gamma(x,y,\lambda, \mu)\\
 & = & \int_{X^2\times [0,1]^2} (\lambda f(x) + \mu g(y)) d\gamma(x,y,\lambda, \mu) \\
 & \leq &   \int_{X^2\times [0,1]^2} \| \lambda x- \mu y\| d\gamma(x,y,\lambda, \mu). \hspace{5cm} \Box
 \end{eqnarray*}
 

 \subsection{The case of probabilities over a simplex} \label{sub23}
 
  We assume here that $X=\Delta(K)$, where $K$ is a non empty finite set. We use $\|p\|=\sum_k |p^k|$ for every vector $p=(p^k)_{k \in K}$ in $\R^K$, and view $X$ as the set of  vectors in $\R^K_+$ with norm 1. 
   $$X=\{p=(p^k)_{k \in K} \in \R^K_+ , \sum_{k \in K} p^k=1\}.$$
 Recall that for $u$ and $v$ in $\Delta(X)$, we have    
 $d_1(u,v) =\sup_{f \in D_1}\;  |u(f)-v(f)|,$ where 
 $ D_1=\{f \in E, \forall x,y \in X, \forall a, b \geq 0, \; a f(x)-b f(y)\leq \| a x-b y\|\}.$

We now introduce an alternative definition of $d_1$ using ``non revealing game functions". These functions come  from the theory of repeated games with incomplete information {\it à la} Aumann Maschler (1995), and the interest  for  the distance  $d_0$  emerged several   years ago while   doing research  on  Markov decision processes with partial observation and repeated games with an informed controller  (see Renault 2011  and  2012a).

 Given a collection of matrices $(G^k)_{k \in K}$ (all of the same finite size $I\times J$) indexed by $K$ and with values in $[-1,1]$, we define the ``non revealing function" $f$ in ${\cal C}(X)$ by:
\begin{eqnarray*}
\forall p\in X, f(p) & =& \Val \left(\sum_{k \in K} p^k G^k\right),\\
 & = & \max_{x \in \Delta(I)}\min_{y \in \Delta(J)} \sum_{i\in I,j\in J}x(i) y(j) \left( \sum_{k \in K} p^k G^k(i,j)\right),\\
 & = & \min_{y \in \Delta(J)}\max_{x \in \Delta(I)} \sum_{i\in I,j\in J}x(i) y(j) \left( \sum_{k \in K} p^k G^k(i,j)\right).
 \end{eqnarray*}
 
 \noindent $f(p)$ is the minmax value of the average matrix $\sum_{k}p^kG^k$. The set of all such non revealing functions $f$, where $I$, $J$ and $(G^k)_{k \in K}$ vary, is denoted by $D_0$. \\
 
 Clearly,  all affine functions from $X$ to $[-1,1]$ belong to $D_0$. 
It is known that the set of non revealing functions is dense in  ${\cal C}(X)$. However, we only consider here non revealing functions defined by matrices with values in $[-1,1]$, and $D_0$ is {\it not} dense in the set of continuous functions from $X$ to $[-1,1]$. As an example, consider  the case where $K=\{1,2\}$ and $f$ in $E$ is piecewise-linear  with $f(1,0)=f(0,1)=0$ and $f(1/2,1/2)=1$. If a function $g$ in $D_0$ is such that $g(1/2,1/2)=1$, then necessarily the values of the two matrix games $G^1$ and $G^2$ are also equal to $1$ since it is the maximum value. Therefore $f$ is not in $D_0$. In fact $f$ is 1-Lipschitz, however  $2 f(1/2,1/2)-f(1,0)=2> \|  2(1/2,1/2)-(1,0)\|=1$, so it is not in $D_1$ which we will see later contains $D_0$ (see lemma \ref{lemdense}). 

 \begin{lem} \label{lemsup}  If  $f$, $g$ belong to $D_0$ and $\lambda\in [0,1]$, then $-f$, $\sup\{f,g\}$, $\inf\{f,g\}$ and $\lambda f + (1-\lambda)g$ are in $D_0$. The linear span of $D_0$ is dense in ${\cal C}(X)$. \end{lem}
 
  \noindent{\bf Proof:}  The proof can be easily deduced from proposition 5.1. page 357 in MSZ, part B. For instance, let $f$ and $g$ in $D_0$ be respectively defined by the collections of matrices $(G^k)_{k \in K}$ with  size $I_1 \times J_1$ and $(H^k)_{k \in K}$ with  size $I_2  \times J_2$. 
  
  Defining  for each $k$, $i_1$, $j_1$: $G'^k(i_1,j_1)= -G^k(j_1,i_1)$  yields a family of matrices $(G'^k)_k$ with size $J_1\times I_1$ inducing $-f$. So $-f\in D_0$.
  
  To get that $\sup\{f,g\}$ belongs to $D_0$, one can  assume w.l.o.g. that $I_1\cap I_2=J_1\cap J_2=\emptyset$. Set $I=I_1\cup I_2$ and $J=J_1\times J_2$. Define for each $k$ the matrix game $L^k$ in $\R^{I \times J}$ by $L^k(i, (j_1,j_2))= G^k(i,j_1)$ if $i\in I_1$, $L^k(i, (j_1,j_2))= H^k(i,j_2)$ if $i\in I_2$. Then for each $p$ in $X$, we have $\Val(\sum_k p^kL^k)= \sup\{f(p),g(p)\}$, so that $\sup\{f,g\}\in D_0$.

\begin{lem} \label{lemdense}  The closure of $D_0$ is $D_1$. 
\end{lem}

  \noindent{\bf Proof:} We first show  that 
  $D_0\subset D_1$. Let $I$ and $J$ be finite sets,  and  $(G^k)_{k \in K}$ be a collection of $I\times J$-matrices with values in $[-1,1]$. Consider  $p$ and $q$ in $X$ and $a$ and $b$ non negative. Then  for all $i$ and $j$:
  $$ \left |\sum_k p^k a G^k(i,j) -\sum_k q^k b G^k(i,j)\right  |   \leq  \sum_k |a p^k -b q^k|= \|a p-bq \|.$$ 
\noindent As a consequence, 
\begin{eqnarray*} 
a \Val \left( \sum_{k \in K} p^k G^k\right) - b \Val\left( \sum_{k \in K} q^k G^k\right) & =  &\Val \left( \sum_{k \in K} a p^k G^k\right) -   \Val\left( \sum_{k \in K} bq^k G^k\right) \\
 & \leq& \|a p-bq \|\end{eqnarray*}

We now show that the closure of $D_0$ is $D_1$. Consider  $f$ in $D_1$, in particular we have $\|f\|_{\infty}\leq 1$. Let $p$ and $q$ be distinct elements in $X$, and define $Y$ as the linear span of $p$ and $q$, and define $\varphi$ from $Y$ to $\R$ such that: $\varphi(\lambda p +\mu  q)=\lambda f(p)+ \mu f(q)$ for all reals $\lambda$ and $\mu$.
 
  If $\lambda \geq 0$ and $\mu \geq 0$, we have $\varphi(\lambda p+ \mu q)\leq \lambda + \mu=\|\lambda p + \mu q\|.$   If $\lambda \geq 0$ and $\mu \leq 0$, we directly use the definition of $D_1$ to get:   $\varphi(\lambda p+ \mu q)\leq \|\lambda p + \mu q\|$. As a consequence, $\varphi$ is a linear form  with norm at most 1 on $Y$. By Hahn-Banach theorem, it can be extended to a linear mapping  on $\R^K$ with the same norm, and we denote by $g$ the restriction of this mapping to $X$. $g$ is affine 
   with $g(p)=\varphi(p)=f(p)$ and $g(q)=\varphi(q)=f(q)$. Moreover, for each $r$ in $X$, we have $\|g(r)\|\leq \|r\|=1$. As a consequence  $g$ belongs to $D_0$.  
   
  Because $D_0$ is stable under the $\sup$ and $\inf$ operations, we can  use    Stone-Weierstrass theorem (see for instance lemma A7.2 in Ash p.392) to conclude that $f$  belongs to the closure  of $D_0$.  \hfill $\Box$
  
  \vspace{0,5cm}

 \begin{defi} Given $u$ and $v$ in $\Delta(X)$, define:  $$d_0(u,v) =\sup_{f \in D_0}\;  u(f)-v(f)$$ \end{defi}

 \begin{pro} \label{pro7} $d_0$ is a distance on $\Delta(X)$ metrizing the weak-* topology. Moreover $d_0=d_1=d_2=d_3$. \end{pro}
 
 \vspace{0,5cm}
 
 \noindent{\bf Proof:}     $d_0=d_1=d_2=d_3$  follows from lemma \ref{lemdense} and theorem \ref{thm1}.  Because the linear  span of  $D_0$ is dense in ${\cal C}(X)$, we obtain the separation property and $d_0$ is a distance on $\Delta(X)$. Because $D_0\subset D_1\subset E_1$, we   have $d_0=d_1\leq d_{KR}$.  Since $(\Delta(X), d_{KR})$ is a compact metric space, the identity map $(\Delta(X), d_{KR})$ to $(\Delta(X), d_0)$ is bicontinuous, and we obtain that $(\Delta(X), d_0)$ is a compact metric space and $d_0$ and $d_{KR}$ are equivalent. (see for instance proposition 2 page 138 Aubin). \hfill $\Box$ 
 
\vspace{0,5cm}
\noindent {\bf Remark:} one can show that allowing for infinite sets $I$, $J$ in the definition of $D_0$ (still assuming that all games $\sum_k p^k G^k$ have  a value) would not change the value of $d_0$.

From now on, we just write $d_*(u,v)$ for the distance $d_0=d_1=d_2=d_3$ on $\Delta(X)$.  Elements of $X$ can be viewed as elements of $\Delta(X)$ (using Dirac measures), and it is well known that for $p$, $q$ in $X$, we have: $d_{KR}(\delta_p, \delta_q)=\|p-q\|.$ We have the same result with $d_*$. 

\begin{lem} For $p$, $q$ in $X$, we have $d_*(\delta_p, \delta_q)=\|p-q\|.$
\end{lem}
\noindent {\bf Proof:} Define $K_1=\{k \in K, p^k \geq q^k\}$, and $K_2=K \backslash K_1$. Consider $f$ affine on $X$ such that $f(k)=+1$ if $k\in K_1$, and $f(k)=-1$ if $k \in K_2$. Then $f \in D_1$, and $d_*(\delta_p, \delta_q)\geq |f(p)-f(q)|=\|p-q\|.$ The other inequality is clear. \hfill $\Box$

\vspace{0,5cm}

We   now present  a dual formulation for our distance, in the spirit of Kantorovich duality formula from optimal transport. For any $u$, $v$ in $\Delta(X)$, we denote by $\Pi(u,v)$ the set of transference plans, or couplings, of $u$ and $v$, that is the set of probability distributions over $X \times X$ with first marginal $u$ and second marginal $v$. Recall (see for instance Villani 2003, p.207):
$$d_{KR}(u,v)=\sup_{f \in E_1}\;  |u(f)-v(f)|= \min_{\gamma \in \Pi(u,v)} \int_{(x,y)\in X \times X} \|x-y\| \; d\gamma(x,y)$$

We will concentrate on probabilities on $X$ with finite support. We denote by $Z=\Delta_f(X)$ the set of such probabilities.

\begin{defi} Let $u$ and $v$ be in $Z$ with respective supports  $U$ and $V$.   We define ${\cal M}_4(u,v)$  as the set $$\left\{(\alpha, \beta)\in ({\R_+}^{U\times V})^2, s.t. \forall x\in U, \forall y \in V, \sum_{y' \in V} \alpha(x,y')= u(x)\; {\rm and}\; \sum_{x' \in U} \beta(x',y) = v(y)\right\}.$$
 $$And \; \;  d_4(u,v)=\inf _{(\alpha, \beta)\in {\cal M}_4(u,v)}  \sum_{(x,y)\in U\times V} \|x\alpha(x,y)-y \beta(x,y)\|  $$
\end{defi}

Notice that diagonal elements in ${\cal M}_4(u,v)$, i.e. measures  $\alpha$ such that $(\alpha, \alpha)\in {\cal M}_4(u,v)$, coincide with elements of $\Pi(u,v)$. 
  ${\cal M}_4(u,v)$ is a  polytope  in  the  Euclidean space $({\R}^{U\times V})^2$, so the infimum in the definition  of $d_4(u,v)$ is achieved. 
  
  \begin{thm} \label{thm2} (Duality formula) Let $u$ and $v$ be in $Z$ with respective supports $U$ and $V$.
   $$d_*(u,v)=  \sup_{f \in D_1}\;  |u(f)-v(f)|= \min_{(\alpha, \beta)\in {\cal M}_4(u,v)}  \sum_{(x,y)\in U\times V} \|x\alpha(x,y)-y \beta(x,y)\|  $$
  where  $D_1=\{f \in E, \forall x,y \in X, \forall a, b \geq 0, \; a f(x)-b f(y)\leq \| a x-b y\|\}$,
  
 \noindent  and ${\cal M}_4(u,v)=\left\{(\alpha, \beta)\in {\R_+}^{U\times V}\times {\R_+}^{U\times V}, s.t. 
\;   \forall (x,y)\in U\times V,  \right.$ 
  
\hspace{2cm}  $\left.\;  \sum_{y' \in V} \alpha(x,y')= u(x)\; {\rm and}\sum_{x' \in U} \beta(x',y) = v(y)\right\}.$
  \end{thm}
    
   The proof is postponed to the next subsection. We conclude this part by a simple but fundamental   property of the distance $d_*$.  
   
   \begin{defi} Given a finite set $S$, we define the {\it posterior} mapping  $\psi_S$ from $\Delta(K\times S)$ to $\Delta(X)$ by:
   $$\psi_S(\pi)=\sum_{s \in S} \pi(s) \delta_{p(s)}$$
   \noindent where for each $s$,  $\pi(s)=\sum_k \pi(k,s)$ and $p(s)=(p^k (s))_{k \in K}\in X$ is the posterior on $K$ given $s$ (defined arbitrarily   if $\pi(s)=0$)
: for each $k$ in $K$, $p^k (s)=\frac{ \pi(k,s)}{\pi(s)}$.    \end{defi}

$\psi_S(\pi)$ is a probability with finite support over $X$. Intuitively, think of  a joint variable  $(k,s)$ being selected according to $\pi$, and an agent just observes $s$. His knowledge on $K$ is then represented by $p(s)$. And $\psi_S(\pi)$ represents the ex-ante information that the agent will know about the variable $k$. 
$\Delta(K \times S)$ is endowed as usual with the $\|.\|_1$ norm. One can show that $\psi_S$ is continuous whenever $X$ is endowed with the weak-* topology.  Intuitively, $\psi_S(\pi)$ has less information than $\pi$, because the agent does not care about $s$ itself but just on the information about $k$ given by $s$.  So one may hope that the mapping $\psi_S$ is 1-Lipschitz (non expansive)  for a well chosen distance on $\Delta(X)$. This is not the case if one uses the Kantorovich-Rubinstein  distance $d_{KR}$, as shown by the example below:

\begin{exa}
Consider the case where $K=\{a,b,c\}$ and $S=\{\alpha,\beta\}$. We denote by $\pi$ and $\pi'$ the following laws on $\Delta(K \times S)$:
\begin{center}
\begin{tabular}{cccc}
 & $S$ & & $S$ \\
$K$ &$\begin{pmatrix}
\frac{1}{4} & 0 \\
0 & \frac{1}{2} \\
\frac{1}{4} & 0
\end{pmatrix}$ & and &$\begin{pmatrix}
\frac{1}{4} & 0 \\
0 & \frac{1}{2} \\
0 & \frac{1}{4} 
\end{pmatrix}$. \\
 & $\pi$ & & $\pi'$
\end{tabular}
\end{center}
Their disintegrations are respectively $\psi_S(\pi)=\frac{1}{2} \begin{pmatrix} \frac{1}{2} \\ 0 \\ \frac{1}{2} \end{pmatrix}+ \frac{1}{2} \begin{pmatrix} 0 \\ 1 \\ 0 \end{pmatrix}$ and $\psi_S(\pi')=\frac{1}{4} \begin{pmatrix} 1 \\ 0 \\ 0 \end{pmatrix}+ \frac{3}{4} \begin{pmatrix} 0 \\ \frac{2}{3}  \\ \frac{1}{3} \end{pmatrix}.$
We define the test function $f: \Delta(K)\rightarrow [-1,1]$ by
\begin{eqnarray*}
f\begin{pmatrix} 0 \\ 1 \\ 0 \end{pmatrix} = \frac{1}{3} &, & f\begin{pmatrix} \frac{1}{2} \\ 0 \\ \frac{1}{2} \end{pmatrix} = -\frac{1}{3}, \\
f\begin{pmatrix} 0 \\ \frac{2}{3} \\ \frac{1}{3} \end{pmatrix} = 1  & and & f\begin{pmatrix} 1 \\ 0 \\ 0 \end{pmatrix} =  \frac{2}{3}.
\end{eqnarray*}
We have $\|\pi-\pi'\|=\frac{1}{2}$ and since $f$ is $1$-Lipschitz, $d_{KR}(\psi_S(\pi),\psi_S(\pi'))\geq \psi_S(\pi')(f)-\psi_S(\pi)(f)=\frac{11}{12}-0> \frac{1}{2}.$ The posterior mapping $\psi_S$ is not $1$-Lipschitz from $(\Delta(K\times S),\|.\|_1)$ to $(\Delta(X),d_{KR})$ .
\end{exa}

However, the next proposition shows that the distance $d_*$ has the desirable property.

\begin{pro}\label{pro3} For each finite set $S$, the mapping $\psi_S$ is 1-Lipschitz   from $(\Delta(K\times S), \|.\|_1)$ to $(\Delta_f(X), d_*)$.

 Moreover, $d_*$ is the largest distance on $Z$  having this property: given $u$ and $v$ in $Z$, we have 
$$d_*(u,v)=\inf \{ \|\pi-\pi'\|_1 , \; s.t. \; {S \; {\rm \it finite}}, \;   \psi_S(\pi)=u,\;  \psi_S(\pi')=v\}.$$
\end{pro}

\noindent {\bf Proof:}  First fix $S$ and $\pi$, $\pi'$ in $\Delta(K \times S)$. Write $u=\psi_S(\pi)$, $u'=\psi_S(\pi')$. For any $f$ in $D_1$, we have:
\begin{eqnarray*}
u(f)-u'(f) &= & \sum_{s\in S} \left (\pi(s) f(p(s))-  \pi'(s) f(p'(s)\right) \\
 & \leq & \sum_{s \in S} \| \pi(s) p(s)-\pi'(s) p'(s)\|\\
 & \leq & \sum_{s\in S} \| (\pi(k,s))_k -(\pi'(k,s))_{k}\| \\
 & \leq &\sum_{s \in S} \sum_{k \in K} | \pi(k,s)-\pi'(k,s)|= \|\pi-\pi'\|_1.
 \end{eqnarray*}
 \noindent So $d_*(u,u')\leq \|\pi-\pi'\|_1$, and $\psi_S$ is 1-Lipschitz. \\
 
Let now $u$ and $v$ be in $Z$. There exists $(\alpha,\beta)\in {\cal M}_4(u,v)$ such that 
\[ d_*(u,v)=\sum_{(x,y)\in U\times V} \|\alpha(x,y)x-\beta(x,y)y\|.\]
Define $S=U\times V$ and $\pi,\pi' \in \Delta(K\times S)$ by
$\pi(k,(x,y))=x(k)\alpha(x,y)$ and $\pi'(k,(x,y))=y(k)\beta(x,y)$. By definition of ${\cal M}_4(u,v)$, $\pi$ and $\pi'$ are probabilities and
\begin{align*}
 \|\pi-\pi'\|_{1,K \times S} & =\sum_{k \in K, (x,y)\in U\times V} |x(k)\alpha(x,y)- y(k)\beta(x,y)|\\
                             & =\sum_{(x,y) \in U\times V} \|\alpha(x,y)x-\beta(x,y)y\| .
\end{align*}  \hfill $\Box$

  \subsection{Proof of the duality formula}

Let $u$ and $v$ be in $\Delta(X)$,  and denote by $U$ and $V$ the respective supports of $u$ and $v$.   We write  $S=X^2\times [0,1]^2$, and we start 
with a lemma, where no finiteness assumption on $U$ or $V$ is needed.

 \begin{lem} \label{lem4}  For each  $\gamma \in {\cal M}_3(u,v)$, we have:
 $$ \int_{ X^2\times [0,1]^2} \|\lambda x-\mu y\| d\gamma(x,y,\lambda, \mu)= 2 + \int_{U \times V \times [0,1]^2} \left( \|\lambda x-\mu y\| - \lambda -\mu  \right) d\gamma(x,y,\lambda, \mu) .$$
  \end{lem}
   \noindent {\bf Proof:} Write $A(\gamma)= \int_{S} \|\lambda x-\mu y\| d\gamma(x,y,\lambda, \mu).$ By definition of $ {\cal M}_3(u,v)$, we have:
   $$\int_{S} \lambda \mathbf{1}_{x \notin U} d\gamma=0, \;{and}\;\int_{S} \mu \mathbf{1}_{y \notin V} d\gamma=0.$$
   So that $\lambda \mathbf{1}_{x \notin U} =\mu \mathbf{1}_{y \notin V}=0$ $\; \gamma$. a.e.
   We can write: $$A(\gamma)= \int_{S} \1_{x\in U, y \in V} \|\lambda x-\mu y\| d\gamma(x,y,\lambda, \mu) +  \int_{S}  \1_{x \in U, y \notin V} \|\lambda x-\mu y\| d\gamma(x,y,\lambda, \mu)$$ $$+  \int_{S} \1_{x \notin U, y \in V}\|\lambda x-\mu y\| d\gamma(x,y,\lambda, \mu)+ \int_{S}\1_{x \notin U, y \notin V} \|\lambda x-\mu y\| d\gamma(x,y,\lambda, \mu)$$
   $$= \int_{S} \1_{x\in U, y \in V} \|\lambda x-\mu y\| d\gamma(x,y,\lambda, \mu) + \int_{S}  \1_{x \in U, y \notin V} \lambda d\gamma(x,y,\lambda, \mu)+  \int_{S} \1_{x \notin U, y \in V} \mu  d\gamma(x,y,\lambda, \mu)+0.$$
   We also have by definition of $ {\cal M}_3(u,v)$ that $1= \int_S \1_{x\in U}\lambda d\gamma$, so that:
   $$1= \int_S \1_{x\in U, y \in V}\lambda d\gamma+ \int_S  \1_{x\in U, y \notin V}\lambda d\gamma.$$
And similarly  $1= \int_S \1_{x\in U, y \in V}\mu d\gamma+ \int_S  \1_{x\notin U, y \in V}\mu d\gamma.$ We obtain:
$$A(\gamma)=2+\int_{S} \1_{x\in U, y \in V} \|\lambda x-\mu y\| d\gamma(x,y,\lambda, \mu)- \int_S  \1_{x\in U, y \in V}\lambda d\gamma- \int_S  \1_{x\in U, y \in V}\mu d\gamma.$$
$\;$ \hfill $\Box$

\vspace{0,5cm}

We assume in the sequel that $U$ and $V$ are finite, and define $d_5(u,v)$  as follows:

\begin{defi}   Define $${\cal M}_5(u,v)=\big\{(\alpha, \beta)=(\alpha(x,y), \beta(x,y))_{(x,y)\in U\times V}\in ({\R}^{U\times V})^2, s.t. \forall x\in U, \forall y \in V,$$
$$ \alpha(x,y)\geq 0, \beta(x,y)\geq 0,   \sum_{y' \in V} \alpha(x,y')\leq u(x)\; {\rm and}\; \sum_{x' \in U} \beta(x',y) \leq v(y)\big\}.$$
 $$And \;   d_5(u,v)=\inf_{(\alpha, \beta)\in {\cal M}_5(u,v)} 2 + \sum_{(x,y)\in U\times V}\left( \|x\alpha(x,y)-y \beta(x,y)\| -\alpha(x,y)-\beta(x,y)\right).$$
 \end{defi} ${\cal M}_5(u,v)$ is a  polytope  in  the  Euclidean space $({\R}^{U\times V})^2$, so the infimum in the definition  of $d_5(u,v)$ is achieved. 
\begin{lem} \label{lem5}$d_3(u,v)\geq d_5(u,v).$
\end{lem}

\noindent{\bf Proof:} Let $\gamma$ be in ${\cal M}_3(u,v)$. Fix for a while $(x,y)$ in $U\times V$, and assume that $\gamma(x,y)>0$. We define $\gamma(.|x,y)$ the conditional probability on $[0,1]^2$ given $(x,y)$ by: for all $\varphi\in C([0,1]^2)$,
$$\int_{[0,1]^2} \varphi(\lambda, \mu) d\gamma(\lambda, \mu|x,y)=\frac{1}{\gamma(x,y)} \int_{(x',y',\lambda, \mu)\in S} \1_{x'=x,y'=y} \varphi(\lambda,\mu) d\gamma(x',y',\lambda, \mu).$$
\indent So that  $$\gamma(x,y) \int_{[0,1]^2}  ( \|\lambda x- \mu y\| -\lambda-\mu)  d\gamma(\lambda, \mu|x,y)=  \int_{(\lambda, \mu)\in [0,1]^2}    (\|\lambda x- \mu y\| -\lambda-\mu) d\gamma(x,y,\lambda, \mu).$$

 The mapping $\Psi: (\lambda, \mu)\mapsto \|\lambda x- \mu y\| -\lambda-\mu$ is convex so by Jensen's inequality we get: $$\int_{(\lambda, \mu)\in [0,1]^2} (\|\lambda x -\mu y\|-\lambda -\mu)  d\gamma(\lambda, \mu|x,y)\geq$$
$$
 \| x \int_{(\lambda, \mu)\in [0,1]^2} \lambda  d\gamma(\lambda, \mu|x,y) - y \int_{(\lambda, \mu)\in [0,1]^2} \mu d\gamma(\lambda, \mu|x,y)\|  $$
  
  $$
  - \int_{(\lambda, \mu)\in [0,1]^2} \lambda  d\gamma(\lambda, \mu|x,y)- \int_{(\lambda, \mu)\in [0,1]^2} \mu  d\gamma(\lambda, \mu|x,y).$$
 We write: $$P(x,y)= \int_{(\lambda, \mu)\in [0,1]^2} \lambda  d\gamma(\lambda, \mu|x,y)\; {and}\; Q(x,y)=\int_{(\lambda, \mu)\in [0,1]^2} \mu  d\gamma(\lambda, \mu|x,y),$$
 so that $$\int_{(\lambda, \mu)\in [0,1]^2} (\|\lambda x -\mu y\|-\lambda -\mu)  d\gamma(\lambda, \mu|x,y)\geq \| xP(x,y) - y Q(x,y)\| -P(x,y)-Q(x,y).  $$
 
 Now, by lemma \ref{lem4} 
 \begin{eqnarray*} A(\gamma) &=& 2 + \sum_{x\in U, y \in V} \int_{(\lambda, \mu)\in [0,1]^2}  \left( \|\lambda x-\mu y\| - \lambda -\mu  \right) d\gamma(x,y,\lambda, \mu)\\
  &=& 2 + \sum_{x\in U, y \in V, \gamma(x,y)>0} \int_{(\lambda, \mu)\in [0,1]^2}  \left( \|\lambda x-\mu y\| - \lambda -\mu  \right) d\gamma(x,y,\lambda, \mu)\\
  & \geq&2 + \sum_{x\in U, y \in V, \gamma(x,y)>0} \gamma(x,y) \left( \| xP(x,y) - y Q(x,y)\| -P(x,y)-Q(x,y)\right).  \end{eqnarray*}
  
For $(x,y)$ in $U\times V$,  define $\alpha(x,y)=\gamma(x,y)P(x,y)\geq 0 \; {and} \; \beta(x,y)=\gamma(x,y)Q(x,y)\geq 0$ 
(with  $\alpha(x,y)=\beta(x,y)=0$ if $\gamma(x,y)=0$). We get:
$$A(\gamma)\geq 2+  \sum_{x\in U, y \in V} \left( \| x\alpha(x,y) - y \beta(x,y)\| -\alpha(x,y)-\beta(x,y)\right).$$
And we have, for each $x$ in $U$:
 \begin{eqnarray*}\sum_{y \in V}\alpha(x,y)&=&\sum_{y \in V, \gamma(x,y)>0} \int_{(\lambda, \mu)\in [0,1]^2} \lambda d\gamma(x,y,\lambda, \mu)\\
  & \leq &  \int_{(y,\lambda, \mu)\in X\times [0,1]^2} \lambda d\gamma(x,y,\lambda, \mu)=u(x).
   \end{eqnarray*}
where the last equality comes from the definition of ${\cal M}_3(u,v)$. Similarly, for each $y$ in $V$ we can show that $\sum_{x \in U} \beta(x,y)\leq v(y)$, and lemma \ref{lem5} is proved. \hfill $\Box$
\begin{lem} $d_5(u,v)\geq d_4(u,v).$ \end{lem}

\noindent{\bf  Proof:}  Consider $(\alpha^*, \beta^*)$ achieving the minimum in the definition of $d_5(u,v)$. Assume that there exists $x^*$ such that $\sum_{y\in V} \alpha(x^*,y)<u(x^*).$ For any $x$ in $X$ and $z$ in $\R^K_+$, one can check that the mapping $l: (\alpha \mapsto \|x\alpha-z\|-\alpha)$ is nonincreasing from $\R_+$ to $\R$ (as the sum of the mappings $l_k: (\alpha \mapsto|\alpha x^k-z^k|-\alpha x^k)$, each  $l^k$ being non increasing in $\alpha$).   As a consequence, one can choose any $y^*$ in $V$ and increase $\alpha(x^*,y^*)$ in order to saturate the constraint  without increasing the objective.  So we can assume without loss of generality that $\sum_{y\in V} \alpha(x^*,y)=u(x^*)$ for all $x^*$ and similarly  $\sum_{x\in U} \beta(x,y^*)=v(y^*)$ for all $y^*$. 

Consequently, 
\begin{eqnarray*}
d_5(u,v) & = &2 + \sum_{(x,y)\in U\times V}\left( \|x\alpha^*(x,y)-y \beta^*(x,y)\| -\alpha^*(x,y)-\beta^*(x,y)\right)\\
 &=& \sum_{(x,y)\in U\times V} \|x\alpha^*(x,y)-y \beta^*(x,y)\| \geq d_4(u,v).
 \end{eqnarray*}
 
\begin{lem} $d_4(u,v)\geq d_2(u,v).$
\end{lem}

\noindent{\bf  Proof:}  Fix $(f,g)\in D_2$ and $(\alpha, \beta)\in  {\cal M}_4(u,v).$
  \begin{eqnarray*}
u(f)+v(g) &=& \sum_{x\in U} f(x)u(x)+ \sum_{y\in Y} g(y) v(y)\\
& = & \sum_{(x,y)\in U \times V} f(x) \alpha(x,y)+g(y)\beta(x,y) \\
 & \leq &\sum_{(x,y)\in U \times V}\|\alpha(x,y)x-\beta(x,y)y\|\leq d_4(u,v).
  \end{eqnarray*}
  
  \vspace{0,5cm}

 We have shown that $d_3(u,v)\geq d_5(u,v)\geq d_4(u,v)\geq d_2(u,v)=d_3(u,v)=d_1(u,v)$. This ends the proof of theorem \ref{thm2}. 


 
 \section{Long-term values for compact non expansive Markov Decision Processes}

 In this section we consider Markov Decision Processes, or Controlled Markov Chains,  with bounded payoffs and transitions with finite support. We will consider   two  closely related models of MDP  and prove   in each case the existence and a characterization for  a general notion of long-term value. The first model deals with  MDP without any explicit action set (hence, payoffs only depend on the current state), such MDP will be called {\it gambling houses} using the terminology of gambling theory (see Maitra and Sudderth 1996).  We will assume in this setup that the set of states $X$  is metric compact and that  the transitions are non expansive with respect to the $KR$-distance on $\Delta(X)$. Since we only use the $KR$-distance here, the theorem for the first model, namely theorem \ref{thm3},  does not use the distance for belief spaces studied in section \ref{sec2}. The second model is the standard model of Markov Decision Processes  with states, actions, transitions  and payoffs, and we will assume that the state space $X$ is a compact subset of a simplex $\Delta(K)$. We will need for this second case an assumption of non expansiveness for the transitions which is closely related to the distance $d_*$ introduced in section \ref{sec2}, see theorem \ref{thm4} later. The   applications in sections \ref{sub42} and \ref{sub43} will be based on the second model.

 \subsection{Long-term  values for  Gambling Houses} \label{sub31}
 
 In this section we consider Markov Decision Processes   of the following form. There is a non empty set of states $X$,  a transition given by a multi-valued mapping  $F:X  \rightrightarrows \Delta_f(X)$ with non empty values, and a payoff (or reward) function $r:X\rightarrow [0,1]$. 
 The idea is that given an initial state  $x_0$ in $X$, a decision-maker (or player)   can choose a probability with finite support $u_1$ in $F(x_0)$, then $x_1$ is selected according to $u_1$ and there is a payoff $r(x_1)$. Then the player has to  select $u_2$ in $F(x_1)$, $x_2$ is selected according to $u_1$ and the player receives the payoff $r(x_2)$, etc... Note that there is no explicit action set here, and that the transitions take values in $\Delta_f(X)$ and hence all have finite support.

We say that  $\Gamma=(X, F, r)$ is a Gambling House.  We assimilate the elements in $X$ with their Dirac measures in $\Delta(X)$, and in case the values of $F$ only consist of  Dirac measures on $X$, we view $F$ as a correspondence from $X$ to $X$ and say that $\Gamma$ is   a {\it deterministic} Gambling House (or a Dynamic Programming problem).  In general we write $Z=\Delta_f(X)$, and an element in $Z$ is written $u=\sum_{x \in X} u(x) \delta_x$. The set of stages   is  $\N^*=\{1,...,t,....\}$, and  a probability distribution over stages is called an evaluation.  Given an evaluation   $\theta=(\theta_t)_{t \geq 1}$   and an initial stage $x_0$ in $X$, the $\theta$-problem $\Gamma_{\theta}(x_0)$  is the  problem induced by a decision-maker starting from $x_0$ and maximizing the expectation of $\sum_{t\geq 1} \theta_t r(x_t)$.

 Formally, we first linearly   extend $r$ and $F$   to $\Delta_f(X)$    by defining for each $u=\sum_{x \in X} u(x) \delta_x$ in $Z$,   the payoff $r(u)=\sum_{x \in X} r(x) u(x)$ and the transition $F(u)=\{\sum_{x \in X} u(x)f(x), s.t. \ f:X \rightarrow Z \text{ and } f(x)\in F(x) \forall x \in X\}$.  We also    define  the {\it mixed extension} of $F$  as  the correspondence  from $Z$  to itself   which associates to every $u=\sum_{x \in X} u(x) \delta_x$ in $\Delta_f(X)$ the  image:  $$\hat{F}(u)=\left\{\sum_{x \in X} u(x) f(x), \; s.t. \ f: X \rightarrow Z \text{ and } f(x)\in \conv F(x)\;  \forall x \in X\right\}.$$

 \noindent  The graph of $\hat{F}$ is the convex hull of the graph of $F$. Moreover $\hat{F}$ is an affine correspondence, as shown by the lemma below.
 
 \begin{lem} $\forall u, u' \in Z$, $\forall \alpha \in [0,1]$,  $\hat{F}(\alpha u +(1-\alpha)u') = \alpha \hat{F}(u)+(1-\alpha) \hat{F}(u').$
  \end{lem}
 \noindent{\bf Proof:} The $\subset $ part is clear. To see the reverse inclusion, let $v= \alpha \sum_{x \in X} u(x) f(x) + (1-\alpha) \sum_{x \in X} u'(x) f'(x)$ be in $\alpha \hat{F}(u)+(1-\alpha) \hat{F}(u')$,  with transparent notations. Define $$g(x)=\frac{\alpha u(x) f(x) + (1 -\alpha) u'(x) f'(x)}{\alpha  u(x)+(1-\alpha) u'(x)},$$ \noindent for each $x$ such that the denominator is positive. Then $g(x)\in \conv F(x)$, and $$v= \sum_{x \in X} (\alpha u(x)+(1-\alpha) u'(x)) g(x)\in \hat{F}(\alpha u+(1-\alpha)u').$$

 \begin{defi} A pure play, or deterministic play,  at $x_0$ is a sequence $\sigma=(u_1,...,u_t,...)\in Z^\infty$ such that  $u_1\in F(x_0)$ and $u_{t+1} \in F(u_{t})$ for each $t\geq 1$. A play, or mixed  play, at $x_0$ is a sequence $\sigma=(u_1,...,u_t,...)\in Z^\infty$ such that $u_1\in \conv F(x_0)$ and  $u_{t+1} \in \hat{F}(u_{t})$ for each $t\geq 1$. We denote by $\Sigma(x_0)$ the set of mixed plays at $x_0$. \end{defi}
 
\noindent  A pure play  is a particular case of a mixed play.   Mixed plays  corresponds to situations where the decision-maker can select, at every stage $t$ and state $x_{t-1}$,  {\it randomly}  the law $u_t$  of the new state. A mixed play at $x_0$ naturally induces a probability distribution over the set ${(X \times \Delta_f(X))}^{\infty}$ of sequences $(x_0,u_0,x_1,u_1,...)$, where $X$ and $Z$ are endowed with the discrete $\sigma$-algebra and ${(X \times \Delta_f(X))}^{\infty}$ is  endowed   with the product $\sigma$-algebra.
 
 \begin{defi}
 Given  an evaluation  $\theta$,  the $\theta$-payoff of a  play $\sigma=(u_1,...,u_t,...)$ is defined as: $\gamma_{\theta}(\sigma)=\sum_{t\geq 1} \theta_t r(u_t)$, and the $\theta$-value at $x_0$ is:
 $$v_{\theta}(x_0)=\sup_{\sigma \in \Sigma(x_0)}  \gamma_{\theta}(\sigma).$$
   \end{defi}
 \noindent It is easy to see that the supremum in the definition of $v_\theta$ can be taken over the set of pure plays at $x_0$. We have the following recursive formula. For each evaluation $\theta=(\theta_t)_{t\geq 1}$ such that $\theta_1<1$, we denote by $\theta^+$ the ``shifted" evaluation $(\frac{\theta_{t+1}}{1- \theta_1})_{t \geq 1}$. We extend linearly $v_{\theta}$ to $Z$, so that the recursive formula can be written: 
 $$\forall \theta \in \Delta(\N^*), \forall x \in X,  \; v_{\theta}(x)= \sup_{u \in \conv F(x)} \left( \theta_1 r(u) + (1-\theta_1) v_{\theta^+} (u)\right).$$
 \noindent And by linearity the supremum can be taken over $F(x)$.   It is also  easy to see that for all evaluation $\theta$ and initial state $x$, we have the inequality:
   \begin{equation} \label{eq1} |v_{\theta}(x) -\sup_{u \in F(x)} v_{\theta} (u)|\leq \theta_1+ \sum_{t\geq 2} | \theta_t- \theta_{t-1}|.\end{equation} 
 
 \vspace{0,5cm}
 
 In this paper, we are interested in the limit behavior  when the decision-maker is very patient. Given an evaluation $\theta$, we define :$$I(\theta)=\sum_{t\geq 1} | \theta_{t+1}-\theta_t|$$
  \noindent  The decision-maker is considered as patient whenever   $I(\theta)$ is small,  so $I(\theta)$ may be seen as the impatience of $\theta$ (see Sorin, 2002 p. 105 and Renault 2012b). When $\theta=(\theta_t)_{t \geq 1}$ is non increasing, then $I(\theta)$ is just $\theta_1$. A classic  example is when   $\theta=\sum_{t=1}^n \frac{1}{n} \delta_t$, the value $v_{\theta}$ is just denoted $v_n$ and the evaluation corresponds to the average payoff from stage 1 to stage $n$. In this case $I(\theta)=1/n \longrightarrow_{n \to \infty} 0.$ We also have $I(\theta)=1/n$ if   $\theta=\sum_{t=m}^{m+n} \frac{1}{n} \delta_t$ for some  non-negative $m$.   Another  example is the case of discounted payoffs, when $\theta=(\lambda {(1-\lambda)}^{t-1})_{t \geq 1}$ for some discount factor $\lambda \in (0,1]$, and in this case   $I(\theta)=\lambda\longrightarrow_{\lambda \to 0} 0.$ 
  
  \begin{defi} \label{defi17} The Gambling House   $\Gamma=(X, F, r)$ has a general limit value $v^*$  if $(v_{\theta})$ uniformly converges to $v^*$ when $I(\theta)$ goes to zero, i.e.: $$\forall \varepsilon >0, \exists \alpha >0, \forall \theta, \;\; \left( \; I(\theta) \leq \alpha\implies \left( \forall x \in X,  |v_{\theta}(x)-v^*(x)|\leq \varepsilon\right)\; \right).$$ \end{defi}

\noindent The existence of the general limit  value implies in particular that  $(v_n)_n$ and $(v_{\lambda})_{\lambda}$  converge to the same limit when $n$ goes  to $+\infty$ and $\lambda$ goes  to $0$. This  is coherent with the result of Lehrer and Sorin (1992), which states that the uniform convergence of $(v_n)_n$ and $(v_{\lambda})_{\lambda}$ are equivalent.

       \vspace{0,5cm}
In the definition of the general limit value, we require all value functions to be close to $v^*$ when the patience is high, but the plays used   may depend on the precise expression  of $\theta$. In the following definition, we require the same play to be simultaneously optimal for all $\theta$ patient enough.

  \begin{defi} \label{defgenunivalue}The Gambling House   $\Gamma=(X, F, r)$ has a general uniform value   if  it has a general limit value $v^*$ and  moreover  for each $\varepsilon>0$ one can find $\alpha>0$ and for each initial  state $x$ a mixed play $\sigma(x)$ at $x$ satisfying:
 $$ \forall \theta, \;\; \left (I(\theta) \leq \alpha \implies \left(\forall x \in X, \gamma_{\theta}(\sigma(x))\geq v^*(x)-\varepsilon \right)\;\right).$$ \end{defi}

Up to now,   the   literature  in repeated games  has  focused on the evaluations  $\theta=\sum_{t=1}^n \frac{1}{n} \delta_t$ and $\theta=(\lambda {(1-\lambda)}^{t-1})_{t \geq 1}$. The standard (Cesàro)-uniform value can be   defined  by restricting the evaluations to be Cesàro means:  for each $\varepsilon>0$ one can find $n_0$ and for each initial  state $x$ a mixed play $\sigma(x)$ at $x$ satisfying:
 $\forall n \geq n_0,  \forall x \in X, \gamma_{n}(\sigma(x))\geq v^*(x)-\varepsilon .$  Recently,  Renault (2011) considered deterministic Gambling Houses and characterized the uniform convergence of the value functions $(v_n)_n$. He also proved the existence of the  standard Cesàro-uniform value under some assumptions,  including the case where the set of states $X$ is metric precompact,  the transitions are non expansive and the payoff function is uniformly continuous.  As  a  corollary, he proved the existence of the uniform value in Partial Observation Markov Decision Processes   with finite set of states (after each stage the decision-maker just observes a stochastic signal more or less correlated to the 
 new state). \\

We now present our main theorem for   Gambling Houses.  Equation (\ref{eq1}) implies that the  general limit value $v^*$ necessarily has to satisfy some rigidity property. The function $v^*$ (or more precisely its linear extension to $Z$) can only be an ``excessive function"  in the terminology of potential theory (Choquet 1956) and gambling houses (Dubins and Savage   1965, Maitra  and Sudderth  1996). 
\begin{defi} An affine   function $w$ defined on $Z$ (or  $\Delta(X)$)   is said to be {\it excessive}  if  for all $x$ in $X$, $ w(x)\geq \sup_{u \in F(x)} w(u)$. \end{defi}  

\begin{exa}
\rm Let us consider the splitting transition   given by $K$ a finite set, $X=\Delta(K)$ and $\forall x\in X, F(x)=\{u \in \Delta(X), \ \sum_{p\in X} u(p)p= x\}$. Then the function $w$ from $Z=\Delta(X)$ to $[0,1]$ is excessive if and only if the restriction of $w$ to $X$ is concave. Moreover given $u,u' \in \Delta(X)$, $u'\in \hat{F}(u)$ if and only $u'$ is the sweeping of $u$ as defined by Choquet (1956): for all continuous concave functions $f$ from $X$ to $[0,1],\ u'(f) \leq u(f)$.
\end{exa}

Assume now that $X$ is a compact metric space and $r$ is continuous. $r$ is naturally extended to an affine continuous  function on $\Delta(X)$ by $r(u)=\int_{p \in X} r(p) du(p)$ for all Borel probabilities on $X$. In the following definition, we consider the closure of the graph of $\hat{F}$ within  the (compact)  set $\Delta(X\times X)$. 
    \begin{defi} \label{invmea} An element $u$ in $\Delta(X)$ is said to be an invariant measure of  the Gambling House  $\Gamma=(X, F, r)$  if $(u,u)\in \cl (Graph \; \hat{F}).$  The set of invariant measures of $\Gamma$ is denoted by $R$, so that: 
  $$ R=\{u \in \Delta(X), (u,u)\in \cl ({Graph \; \hat{F}})\}.$$ \end{defi}
  \noindent $R$ is a convex compact subset of $\Delta(X)$.  Recall that for $u$ and $u'$ in $\Delta(X)$, the Kantorovich-Rubinstein  distance between $u$ and $u'$ is denoted by $  d_{KR}(u,u')=\sup_{f \in E_1}\;  |u(f)-u'f)|.$

  

\begin{thm} \label{thm3} Consider a Gambling House    $\Gamma= (X, F, r)$ such that $X$ is a  compact metric space, $r$ is continuous and $F$ is non expansive with respect to the KR distance:
$$\forall x \in X, \forall x' \in X, \forall u \in F(x), \exists u'\in F(x') s.t.\; d_{KR}(u,u')\leq d(x,x').$$
Then the Gambling House  has a general   uniform value $v^*$   characterized by:  
\begin{eqnarray*} \forall x \in X, \; v^*(x) & = \inf &\big\{ 
w(x),  w:\Delta(X)\rightarrow [0,1] \;  {\rm affine}\;  C^0 \; s.t. \\
 & \; &  (1)\;   \forall y\in X, w(y)\geq \sup_{u \in F(y)} w(u)  \; {\rm and }\;      (2)  \forall u \in R, w(u)\geq r(u) \;  \big \}.\end{eqnarray*}
That is, $v^*$ is the smallest continuous affine function on $X$ which is  1) excessive and 2) above the running payoff $r$ on invariant measures.\end{thm}    

Notice that: 

1) when  $\Gamma= (X, F, r)$  is deterministic, the hypotheses are satisfied  as soon as $X$ is metric compact for some metric $d$, $r$ is continuous and $F$ is non expansive  for $d$. 

2) when $X$ is finite, one can use the distance $d(x,x')=2$ for all $x\neq x'$ in $X$, so that for $u$ and $u'$ in $\Delta(X)$, $d_{KR}(u,u')=\|u-u' \|_1=\sum_{x \in X} |u(x)-u'(x)|$, and the hypotheses are automatically satisfied. We will prove  later a more general  result for  a model of MDP with finite state space,  allowing for explicit actions influencing transitions and payoffs (see corollary \ref{cor11}).

\begin{rem}
\rm The formula also holds when there is no decision maker, i.e. when $F$ is single-valued, and there are some similarities with the Von Neumann ergodic theorem (1932). Let $Z$ be a Hilbert space and $Q$ be a linear isometry on $Z$, this  theorem states that for all $z\in Z$, the sequence $z_n=\frac{1}{n} \sum_{t=1}^n Q^t(z)$ converges to the projection $z^*$ of $z$ on   the set $R$ of fixed points of $Q$. Using the linearity and the non expansiveness leads to a characterization by the set of fixed points. In particular, having in mind  linear payoff functions of the form  $(z\mapsto <l,z>$), we have that the projection $z^*$ of $z$ on $R$ is  characterized by: 
$$\forall l\in Z, <l,z^*>= <l^*,z>=\inf\{<l',z>, l'\in R \; {\rm and}\; <l',r>\geq <l,r> \forall r \in R\}.$$\end{rem}

\begin{exa} \label{ex39} \rm  We consider here a basic periodic   sequence of 0 and 1. Let $X=\{0,1\}$ and for all $x \in X$, $F(x)=\{1-x\}$ and $r(x)=x$. There is a unique invariant measure $u=1/2 \delta_0 + 1/2 \delta_1$, and  the general uniform value exists and satisifes $v^*(x)=\frac{1}{2}$ for all states $x$. Notice that considering evaluations $\theta=(\theta_t)_t$ such that $\theta_t$ is small for each $t$ without  requiring $I(\theta)$ small,  would not necessarily lead to $v^*$. Consider for instance $\theta^n=\sum_{t=1}^n \frac{1}{n}\delta_{2t}$ for each $n$, we have $v_{\theta^n}(x)=x$ for all $x$ in $X$. 
\end{exa}
\begin{exa}\label{cercle}
\rm The state space is the unit circle, let $X=\{x\in \C, |x|=1\}$ and $F(e^{i\alpha})=e^{i(\alpha+1)}$ for all real $\alpha$. If we denote by $\mu$ the uniform distribution (Haar probability measure) on the circle, the mapping $F$ is $\mu$-ergodic and $\mu$ is $F$-invariant. By Birkhoff's theorem (1931), we know that the time average converges to the space average $\mu$-almost surely. Here $\mu$ is the  unique invariant measure, and we  obtain that the general uniform value is the constant: 
\[\forall x \in X, \; v^*(x)=\frac{1}{2\pi}\int_0^{2\pi} r(e^{i\alpha})d\alpha.\]
\noindent Notice that the value $v_{\theta}(x)$ converges to $v^*(x)$ for all $x$ in $X$, and not only for $\mu$-almost all $x$ in $X$. \end{exa}

\begin{exa} \rm  Let $\Gamma=(X,F,r)$ be a MDP satisfying the hypotheses of the theorem \ref{thm3} such that for all $x\in X$, $\delta_x\in F(x).$ Therefore the set $R$ is equal to $\Delta(X)$. In the terminology of Gambling Theory (see Maitra Sudderth, 1996), $\Gamma$ is called   a {\it leavable}  gambling house since at each stage the player can stay at the current state. The limit value $v^*$ is here  characterized by: 
$$v^*    = \inf  \{ v:X\rightarrow [0,1] \; C^0, v  \; {is} \; {excessive} \; \; and \;{v \geq r}  \}.$$
In the above formula, $v$ excessive means: $\forall x \in X, v(x) \geq \sup_{u \in F(x)} \E_u (v).$ This is a variant  of the {\it fundamental theorem of gambling theory}   (see section 3.1 in Maitra Sudderth 1996).  \end{exa}

\begin{exa} \rm \label{infini} The following deterministic  Gambling House, which is an extension   of example 1.4.4. in Sorin (2002)   and  of  example 5.2 of Renault (2011), shows that the assumptions of theorem \ref{thm3}  allow for many  speeds of convergence to the limit value $v^*$. Here   $l> 1$ is  a fixed parameter,   $X$ is the simplex $\{x=(p^a,p^b,p^c)\in \R^3_+, p^a+p^b+p^ c=1\}$ and the initial state is $x_0=(1,0,0)$. The payoff is $r(p^a,p^b,p ^c)=  p^b-p^c$, and the transition is defined by: $F(p^a,p^b,p ^c)=\{( (1-\alpha-\alpha^l)p^a, p^b+\alpha p^a, p^c+\alpha^lp^a), \alpha \in [0,1/2]\}$.   

The probabilistic interpretation  is the following: there are 3 points $a$, $b$ and $c$, and the initial point is $a$. The payoff is 0 at $a$, it is +1 at $b$, and -1 at $c$. At point $a$, the decision maker has to choose $\alpha\in [0,1/2]$~: then  $b$ is reached  with probability $\alpha$, $c$ is reached with probability $\alpha^l$, and the play  stays in $a$ with the remaining probability $1-\alpha-\alpha^l$. When $b$ (resp. $c$) is reached, the play stays at $b$ (resp. $c$) forever. So the decision maker starting at point $a$ wants to reach $b$ and to avoid $c$. By playing at each stage $\alpha>0$ small enough, he can get as close to $b$ as he wants. 

Back to our deterministic setup, we use norm $\|.\|_1$ and obtain that $X$ is compact, $F$ is non expansive and $r$ is continuous, so that theorem \ref{thm3} applies. The limit value is given by  $v^*(p^a,p^b,p^c)=p^a+p^b$, and if we denote by $x_\lambda$ the value $v_\lambda(x_0)$, we have for all $\lambda \in (0,1]$:
$x_{\lambda}=\phi(x_{\lambda})$,  where for all $x \in \R$, 
$$\phi(x)=\max_{\alpha \in [0,1/2]} (1-\lambda)(1-\alpha-\alpha^l) x+\alpha.$$
Since $x_\lambda \in (0,1)$, the first order condition gives $(1-\lambda) x_{\lambda}(-1- l \alpha^{l-1})+1 = 0$  and we can obtain: 
$$ x_\lambda=\frac{1}{(1-\lambda)}\left(l \left(\frac{\lambda}{(1-\lambda)(l-1)}\right)^{\frac{l-1}{l}}+1 \right)^{-1}.$$

Finally we can compute  an  equivalent of $x_\lambda$ as $\lambda$ goes  to $0$. We have
\[
\left(\frac{\lambda}{(1-\lambda)(l-1)}\right)^{\frac{l-1}{l}}=(\frac{1}{l-1})^{\frac{l-1}{l}}\lambda^{\frac{l-1}{l}}(1+ o(\lambda^{\frac{l-1}{l}}))
\]
so that 
\begin{align*}
v_{\lambda}(x_0)& =(1-\lambda)\frac{1}{l \left( (\frac{1}{l-1})^{\frac{l-1}{l}}\lambda^{\frac{l-1}{l}}+ o(\lambda^{\frac{2l-2}{l}}) \right)+1} \\
v_\lambda(x_0) & = 1- C \lambda^{\frac{l-1}{l}}+ o(\lambda^{\frac{l-1}{l}}) \text{ with }C=\frac{l}{(l-1)^\frac{(l-1)}{l}}.
\end{align*}

\end{exa}
 
\subsection{Long-term  values for  standard MDPs}

  A standard Markov Decision Problem $\Psi$ is given by a non empty set of states $X$, a non empty set of actions $A$, a mapping $q:X\times A\rightarrow \Delta_f(X)$ and a payoff function $g:X\times A\rightarrow[0,1].$ At each stage, the player learns the current state $x$ and chooses an action $a$. He then receives the payoff $g(k,a)$,  a new state is drawn accordingly to $q(k,a)$ and the game proceeds to the next stage. 
  
 \begin{defi} A pure, or deterministic,  strategy   is a sequence of mappings $\sigma=(\sigma_t)_{t \geq 1}$ where $\sigma_t:(X\times A)^{t-1}\rightarrow A$ for each $t$. A strategy (or behavioral strategy) is a sequence of mappings $\sigma=(\sigma_t)_{t \geq 1}$ where $\sigma_t:(X\times A)^{t-1}\rightarrow \Delta_f(A)$ for each $t$.  We denote by  $\Sigma$ the set of strategies.    \end{defi}

A pure strategy is a particular case of strategy. An initial state $x_1$ in $X$ and a strategy $\sigma$ naturally induce a probability distribution with finite support over the set of finite histories  $(X\times A)^{n}$  for all  $n$, which can be uniquely extended to a probability over the set  $(X\times A)^{\infty}$ of infinite  histories.
\begin{defi}
Given an evaluation  $\theta$ and an initial state  $x_1$ in $X$, the $\theta$-payoff of a strategy $\sigma$ at $x_1$  is defined as  $\gamma_{\theta}(x_1, \sigma)=\E_{x_1,\sigma} \left( \sum_{t\geq 1} \theta_t g(x_{t},a_{t}) \right)$, and the $\theta$-value at $x_1$ is:
 $$v_{\theta}(x_1)=\sup_{\sigma \in \Sigma}  \gamma_{\theta}(x_1, \sigma).$$
\end{defi}
As for gambling houses, it is easy  to see that the supremum can be taken over the smaller set of pure strategies, and one can derive a recursive formula linking the value functions.  General limit and uniform values are defined as in the previous subsection \ref{sub31}.

 \begin{defi} \label{defi17}  Let $\Psi=(X, A, q, g)$  be a standard MDP.
 
 $\Psi$  has a general limit value $v^*$  if $(v_{\theta})$ uniformly converges to $v^*$ when $I(\theta)$ goes to zero, i.e. for each $\varepsilon>0$ one can find $\alpha>0$ such that:  $$\forall \theta, \;\; \left( \; I(\theta) \leq \alpha\implies \left( \forall x \in X,  |v_{\theta}(x)-v^*(x)|\leq \varepsilon\right)\; \right).$$  

  $\Psi$ has a general uniform value   if  it has a general limit value $v^*$ and moreover   for each $\varepsilon>0$ one can find $\alpha>0$ and  a behavior strategy $\sigma(x)$  for each initial  state $x$ satisfying:
 $$ \forall \theta, \;\; \left (I(\theta) \leq \alpha \implies \left(\forall x \in X, \gamma_{\theta}(x, \sigma(x))\geq v^*(x)-\varepsilon \right)\;\right).$$ \end{defi}
 
We now present a notion of invariance for the MDP $\Psi$. The next definition will be similar to   definition   \ref{invmea}, however  one needs to be  slightly more sophisticated here to incorporate the payoff component.  Assume now that $X$ is a compact metric space, and define for each $(u,y)$ in $\Delta_f(X)\times [0,1]$, 
 $$\hat{F}(u,y)=\left\{\left(\sum_{x \in X} u(x) q(x, a(x)),\sum_{x \in X} u(x) g(x, a(x))\right), \; where\;  a: X\rightarrow \Delta_f(A) \right\}.$$
 \noindent where $q(x,.)$ and $g(x,.)$ have been linearly extended for all $x$. We have defined a correspondence $\hat{F}$  from $\Delta_f(X)\times [0,1]$ to itself. It is easy to see that $\hat{F}$ always is an affine correspondence (see   lemma \ref{lem*} later).  In the following definition we consider the  closure of the graph of $\hat{F}$ within the compact set $\left(\Delta(X)\times [0,1]\right)^2$, with the weak topology.
 
  \begin{defi} An element $(u,y)$ in $\Delta(X)\times[0,1]$ is said to be an invariant  couple  for the MDP  $\Psi$  if $((u,y),(u,y))\in  cl(Graph(\hat{F})) .$ The set of invariant couples of $\Psi$ is denoted by $RR$. \end{defi}

Our main result for standard MDPs is the following, where $X$ is assumed to be a  compact subset of a simplex  $\Delta(K)$.  Recall that  $D_1=\{f \in {\cal C}(\Delta(K)), \forall x,y \in \Delta(K), \forall a, b \geq 0, \; a f(x)-b f(y)\leq \| a x-b y\|_1\}$, and any   $f$ in $D_1$ is linearly extended to $\Delta(\Delta(K))$. 
\begin{thm}\label{thm4}  
Let   $\Psi=(X,A,q,g)$ be a standard MDP where $X$ is a compact subset of a simplex  $\Delta(K)$, such that: 
$$\forall x \in X, \forall y\in X, \forall a \in A,   \forall f \in D_1,  \forall \alpha \geq 0, \forall \beta \geq 0,$$
$$ |\alpha f(q(x,a))-\beta f(q(y,a))|\leq \|\alpha x -\beta y\|_1 \; {\rm and}\; 
  |\alpha g(x,a)-\beta g(y,a)| \leq \|\alpha x -\beta y\|_1
.$$
then $\Psi$  has  a general uniform value $v^*$ characterized by: for all $x$ in $X$, 
\begin{eqnarray*} v^*(x) & = \inf &\big\{ 
w(x),  w:\Delta(X)\rightarrow [0,1] \;  {\rm affine}\;  C^0 \; s.t. \\
 & \; &  (1)\;   \forall x'\in X, w(x')\geq \sup_{a \in A} w(q(x',a))  \; {\rm and }\;      (2) \;  \forall (u,y) \in RR, w(u)\geq y \;  \big \}.\end{eqnarray*}
 \end{thm}

The proof of theorem  \ref{thm4} will be in section \ref{sec35}.   An immediate corollary is when the state space  is finite. 
 \begin{cor} \label{cor11}
Consider a standard MDP  $(K,A,q,g)$ with a finite set of states $K$. Then it has  a general uniform value $v^*$, and for each state $k$:
\begin{eqnarray*} v^*(k) & = \inf &\big\{ 
w(k),  w:\Delta(K)\rightarrow [0,1] \;  {\rm affine}\;   \; s.t. \\
 & \; &  (1)\;   \forall k'\in K, w(k')\geq \sup_{a \in A} w(q(k',a))  \; {\rm and }\;      (2)  \forall (p,y) \in RR, w(p)\geq y \;  \big \}.\end{eqnarray*}
with 
$RR=\{(p,y)\in \Delta(K) \times [0,1], ((p,y),(p,y)) \in \cl(conv(Graph(F))) \}$ and $F(k,y)=\{(q(k,a),g(k,a)) , a \in A \}.$
\end{cor}

\noindent{\bf Proof:} $K$ is viewed as a subset of the simplex $\Delta(K)$, endowed with the $L^1$-norm. Fix $k$, $k'$ in $K$, $a$ in $A$, $\alpha\geq 0$ and $\beta\geq 0$. We have 
\[
\|\alpha k- \beta k'\|=\begin{cases} |\alpha-\beta| \text{ if } k=k' , \\  \alpha+ \beta \text{ otherwise. }
                                                                       \end{cases}\]

First, 
\[
|\alpha g(k,a)-\beta g(k',a)|\leq \begin{cases} |\alpha-\beta| g(k,a) \text{ if } k=k'   \\  \alpha+ \beta \text{ otherwise}
                                                                       \end{cases}
\]
                                                                       , so   in all cases  $|\alpha g(k,a)-\beta g(k',a)|\leq  \|\alpha k -\beta k'\|.$
Secondly,  consider   $f\in D_1$. $f$ takes values in $[-1,1]$, so similarly   we have: $
|\alpha f(q(k,a))-\beta f(q(k',a))|\leq  \|\alpha k -\beta k'\|$. So we can apply theorem \ref{thm4}, and the graph of $\hat{F}$ is the convex hull of the graph of $F$. $\hfill \Box$

\begin{rem}
\rm When the set of actions is finite, we are in the setting of Blackwell (1962) and the value is characterized by the Average Cost Optimality Equation. In fact in this setting, our characterization leads to a dual formulation of a result of Denardo and Fox (1968). Denardo and Fox (1968) showed that the value $v^*$  is the smallest (pointwise) excessive function for which there exists a vector  $h\in \R^K$ such that $(v^*,h)$ is superharmonic in the sense of Hordjik and Kallenberg (1979) , i.e.
\begin{align}
\forall k \in K, \ a \in A \ v^*(k)+h(k)\geq  r(k,a)+ \sum_{k'} q(k,a)(k') h(k').
\end{align}
Given a function $w$ the existence of a vector $h$ such that $(w,h)$ is superharmonic is a linear programming problem with $K\times A$ inequalities. By Farkas' lemma it has a solution if and only if a dual problem has no solution, and the dual programming problem is to find a solution $\pi \in \R^{K\times A}$  of the following system:
\begin{center}
\begin{tabular}{rrl}
$\forall k\in K$ & $\sum_{a'\in A} \pi(k,a')$ & $= \sum_{k' \in K, a'\in A} \pi(k',a') q(k',a')(k)$ \\
$\forall (k,a)\in K\times A$ & $\pi(k,a) $ & $\geq 0$ \\
$\forall k\in K$ & $\sum_{a'\in A} \pi(k,a) g(k,a')$ & $> v(k)$. \\
\end{tabular}
\end{center}
If we denote by $p$ the marginal of $\pi$ on $K$ and define for all  $k$ such that $p(k)>0$, $\sigma(k)=\frac{\pi(k,a)}{p(k)}$ and set $\sigma(k)$ to any probability otherwise, then $\sigma$ is a strategy in the MDP. Moreover $p$ is invariant under $\sigma$ and the stage  payoff $y$ is greater than $v(p)$, thus the couple $(p,y)$ is in $RR$ and the condition $(2)$ in corollary \ref{cor11} is not satisfied. Reciprocally since the action state is compact, given $(p,y) \in RR$, there exists a strategy $\sigma$ such that $p$ is invariant under $\sigma$ and the payoff is $y$. Therefore if the condition (2) is not true then there exists $h\in \R^k$ such that $(w,h)$ is superharmonic.  Note that Denardo and Fox state a dual of the minimization problem  and obtain an explicit dual maximization problem whose solution is the value. Hordjik and Kallenberg exhibit from the solutions of this dual problem an optimal strategy.
\end{rem}

 
\subsection{Proof    of theorem \ref{thm3}}

In this section we consider a compact metric space $(X,d)$, and we use the Kantorovich-Rubinstein distance $d=d_{KR}$ on $\Delta(X)$.   We write $Z=\Delta_f(X)$, $\overline{Z}=\Delta(X)$. We start with a lemma. 

 \begin{lem} \label{lemKR} Let  $F:X  \rightrightarrows \Delta_f(X)$ be  non expansive for  $d_{KR}$.  Then the mixed extension of $F$ is 1-Lipschitz from $\Delta_f(X)$ to $\Delta_f(X)$ for $d_{KR}$. 
   \end{lem}

   \noindent{\bf Proof of lemma \ref{lemKR}.}   We first show  that the mapping $(p \mapsto\conv F(p))$ is non expansive from $X$ to $Z$. Indeed, consider $p$ and $p'$ in $X$, and $u=\sum_{i \in I} \alpha_i u_i$, with $I$ finite, $\alpha_i \geq 0$, $u_i\in F(p)$ for each $i$, and $\sum_{i \in I}\alpha_i=1$. By assumption for each $i$ one can find $u'_i$ in $F(p')$ such that $d_{KR}(u_i,u'_i)\leq d(p,p').$ Define $u'=\sum_{i \in I} \alpha_i u'_i$ in $\conv F(p')$. We have:
   \begin{eqnarray*}
  d_{KR}(u,u') & = &\sup_{f \in E_1} \left(\sum_{i} \alpha_i u_i(f)- \sum_i \alpha_i u'_i(f)\right), \\
   &= &\sup_{f \in E_1} \sum_{i \in I} \alpha_i (u_i(f)-u'_i(f)),\\
   & \leq & \sum_{i \in I} \alpha_i \; d_{KR}(u_i,u'_i),\\
   & \leq & d(p,p').
   \end{eqnarray*}

    We now prove that $\hat{F}$ is 1-Lipschitz from $Z$ to $Z$. Let $u_1$, $u_2$ be  in $Z$ and $v_1=\sum_{p \in X} u_1(p) f_1(p)$, where $f_1(p)\in \conv F(p)$ for each $p$. By the Kantorovich   duality formula, there exists a coupling $\gamma=(\gamma(p,q))_{(p,q)\in X \times X}$ in $\Delta_f(X \times X)$ with first marginal $u_1$ and second marginal $u_2$ satisfying: $$d_{KR}(u_1,u_2)=\sum_{(p,q)\in X \times X} \gamma(p,q) d(p,q).$$ 
   \noindent For each $p$, $q$ in $X$ by the first part of this  proof there exists $f^p(q) \in \conv F(q)$ such that $d_{KR}(f^p(q), f_1(p))\leq d(p,q)$. We define: 
     $$f_2(q)=\sum_{p \in X} \frac{\gamma(p,q)}{ u_2(q)} f^p(q) \in \conv F(q), \; {\rm and} \;  v_2=   \sum_{q \in X} u_2(q) f_2(q)\in \hat{F}(u_2).$$
   We now conclude. 
    \begin{eqnarray*}
  d_{KR}(v_1, v_2) & = &d_{KR}\left( \sum_{p \in X} u_1(p) f_1(p), \sum_{q \in X} u_2(q) f_2(q)\right) \\
   & = & d_{KR}\left( \sum_{p,q} \gamma(p,q) f_1(p), \sum_{q, p} \gamma(p,q) f^p(q) \right)\\
   & \leq & \sum_{p,q} \gamma(p,q) d_{KR}(f_1(p),f^p(q))\\
   & \leq &  \sum_{p,q} \gamma(p,q) d(p,q)= d_{KR}(u_1,u_2). 
   \end{eqnarray*} 
The mixed extension of $F$ is $1$-lipschitz. $\hfill \Box$

\vspace{1cm}

We now consider a Gambling House $\Gamma=(X, F,r)$ and  assume the hypotheses of theorem \ref{thm3} are satisfied. We  will work\footnote{A   variant of the proof would be to consider the Gambling House  on $\Delta(X)$  where the transition correspondence is defined so that  its graph is the closure of the graph of $\hat{F}$. Part 1) of lemma \ref{lem43} shows this correspondence   is also non expansive.} with  the deterministic Gambling House  $\hat{\Gamma}=(\Delta_f(X),\hat{F},r)$.    Recall that $r$ is extended to an affine and continuous mapping  on $\Delta(X)$ whereas $\hat{F}$ is an affine non expansive correspondence from $Z$ to $Z$. \\

 For $p$ in $X$, the  pure   plays in $\hat{\Gamma}$ at the  initial state $\delta_p$ coincide with the mixed plays in ${\Gamma}$ at the  initial state $p$. As a consequence,  the $\theta$-value for $\Gamma$ at $p$  coincides with the $\theta$-value for $\hat{\Gamma}$ at $\delta_p$,  which  is written  $v_{\theta}(p)= v_{\theta}(\delta_p)$. Because  $\hat{F}$ and $r$ are affine on $Z$,   the $\theta$-value  for $\hat{\Gamma}$, as a function defined on $Z$, is the affine extension of the original $v_\theta$ defined on  $X$. So we have a unique value function $v_{\theta}$ which is defined on $Z$ and is affine. Because $\hat{F}$ is 1-Lipschitz and $r$  is uniformly continuous, all the value functions $v_{\theta}$ have the same modulus of continuity as $r$, so $(v_{\theta})_{\theta}$ is an equicontinuous family of mappings from $Z$ to $[0,1]$. Consequently, we  extend $v_{\theta}$ to an affine mapping on $\overline{Z}$ with the same modulous of continuity, and the family $(v_{\theta})_{\theta}$ now is an equicontinuous\footnote{$Z$ being precompact, this is enough to obtain the existence of a general limit value, see Renault 2012b. Here we will moreover obtain a characterization of this value and the existence of the general uniform value.}
 family of mappings from $\overline{Z}$ to $[0,1]$.

We define $R$ and $v^*$ as in the statements of theorem \ref{thm3}, so that  for all $x$ in $X$,
\begin{eqnarray*} v^*(x) & = \inf &\big\{ 
w(x),  w:\overline{Z}\rightarrow [0,1] \;  {\rm affine}\;  C^0 \; s.t. \\
 & \; &  (1)\;   \forall y\in X, w(y)\geq \sup_{u \in F(y)} w(u)  \; {\rm and }\;      (2)  \forall u \in R, w(u)\geq r(u) \;  \big \}.\end{eqnarray*}
We start with a technical lemma using the non-expansiveness of $\hat{F}$. 
  
  \begin{lem} \label{lem43} 1) Given $(u,u')$ in $\cl (Graph(\hat{F}))$, $v$ in $Z$ and $\varepsilon>0$, there exists $v'\in \hat{F}(v)$ such that $d(u',v')\leq d(u,v)+ \varepsilon$.
  
%
2)   Given a sequence $(z_t)_{t\geq 0}$ of elements of $\overline{Z}$ such that $(z_{t},z_{t+1})\in \cl (Graph(\hat{F}))$ for all $t\geq 1$, for each $\varepsilon$ one can find a sequence $(z'_t)_{t\geq 0}$ of elements of $Z$ such that $(z'_t)_{t \geq 1}$ is a play at $z'_0$,  and $d(z_t,z'_t)\leq \varepsilon$ for each $t\geq 0$. 
  \end{lem}
  
  \noindent{\bf Proof of lemma \ref{lem43}:} 1) For all $\varepsilon>0$ there exists $(z,z')\in Graph(\hat{F})$ such that $d(z,u)\leq \varepsilon$ and $d(z',u')\leq \varepsilon$. Because $\hat{F}$ is non expansive, one can find $v'$ in $\hat{F}(v)$ such that $d(z',v')\leq d(z,v)$. Consequently, $d(v',u')\leq d(v',z')+d(z',u')\leq d(z,v) +\varepsilon\leq d(u,v)+2 \varepsilon$.  
  
%

  2) It is first easy to construct  $(z'_0,z'_1)$ in the graph of $\hat{F}$ such that $d(z'_0, z_0)\leq \varepsilon$ and $d(z'_1, z_1)\leq \varepsilon$.  $(z_1,z_2)\in  \cl (Graph(\hat{F}))$ so by 1) one can find $(z'_2)$ in $\hat{F}(z'_1)$ such that $d(z_2, z'_2)\leq d(z_1,z'_1)+ \varepsilon ^2\leq \varepsilon+ \varepsilon ^2$.  Iterating, we  construct a play $(z'_t)_{t\geq 1}$ at $z'_0$ such that $d(z_t,z'_t)\leq \varepsilon + \varepsilon^2+ ...+ \varepsilon^t$ for each $t$.

    \begin{pro} \label{pro4} $\Gamma$ has a general limit value given by $v^*$. \end{pro}

\noindent{\bf Proof of proposition \ref{pro4}:}   By Ascoli's theorem, it is enough  to show that any limit point of $(v_{\theta})_{\theta}$ (for the  uniform convergence) coincides with $v^*$.  We thus assume that $(v_{\theta^k})_k$ uniformly converges to $v$ on $\overline{Z}$ when $k$ goes to $\infty$, for a family of evaluations satisfying: $$\sum_{t\geq 1} |\theta^k_{t+1}-\theta^k_t|\longrightarrow_{k \to \infty} 0.$$

\noindent And we  need to show that $v=v^*$. \\

A) We first show that $v\geq v^*$.

 It is plain that $v$ can be extended to an  affine function on $\overline{Z}$ and has the same modulus of continuity of $r$. Because $\sum_{t\geq 1} |\theta^k_{t+1}-\theta^k_t|\longrightarrow_{k \to \infty} 0$, we have by equation (\ref{eq1}) of section \ref{sub31} that:  $ \forall y\in X, v(y)= \sup_{u \in F(y)} v(u)$. 

  Let now $u$ be in $R$. By lemma \ref{lem43} for each $\varepsilon$  one can find $u_0$ in $Z$ and a play $(u_1,u_2,...,u_t,...)$ such that   $u_t\in \hat{F}(u_{t-1})$ and $d(u,u_t)\leq \varepsilon$ for all $t\geq 0$. Because $r$ is uniformly continuous, we get  $v(u)\geq r(u)$. 
  
 By definition of $v^*$ as an infimum, we obtain: $v^*\leq v$.\\
 
 B) We show that $v^*\geq v$. Let $w$ be a continous affine mapping from $\overline{Z}$ to $[0,1]$ satisfying $(1)$ and $(2)$ of the definition of $v^*$.  It is enough to show that $w(p)\geq v(p)$ for each $p$ in $X$.  Fix $p$ in $X$ and $\varepsilon>0$.
 
 For each $k$, let $\sigma^k=(u^k_1,...,u^k_t,...)\in Z^{\infty}$ be a play at $\delta_p$ for $\hat{\Gamma}$ which is almost optimal for the  $\theta^k$-value, in the sense that $\sum_{t \geq 1} \theta^k_t r(u^k_t) \geq v_{\theta^k}(p)-\varepsilon$. Define:$$u(k)=\sum_{t=1}^{\infty} \theta^k_t u^k_t \in \overline{Z},\;{\rm and }\;  u'(k)=\sum_{t=1}^{\infty} \theta^k_t u^k_{t+1} \in \overline{Z}.$$
 \noindent $u(k)$ and $u'(k)$ are well-defined limits of  normal  convergent series in the Banach space ${{\cal C}(X)}'$. Because $\hat{F}$ is affine, its graph is a convex set and $(u(k),u'(k))\in \cl(Graph(\hat{F}))$ for each $k$. 
 
 Moreover,  we have $d(u(k),u'(k))\leq {\rm diam}(X) (\theta_1^k + \sum_{t=2}^\infty |\theta^k_t-\theta^k_{t-1}|)$, where ${\rm diam}(X)$ is the diameter of $X$. 
 Consequently,  $\sum_{t\geq 1} |\theta^k_{t+1}-\theta^k_t|\longrightarrow_{k \to \infty} 0$ implies $d(u(k),u'(k))\longrightarrow_{k \to \infty} 0$. Considering a limit point of the sequence $(u(k),u'(k))_k$, we obtain some $u$ in $R$. By assumption on $w$, $w(u)\geq r(u)$. Moreover, for each $k$ we have $r(u(k))=\sum_{t \geq 1} \theta^k_t r(u^k_t) \geq v_{\theta^k}(p)-\varepsilon$, so $r(u)\geq v(p)-\varepsilon$.
 
Because $w$ is excessive, we obtain that for each $k$ the sequence $(w(u^k_t))_t$ is non increasing, so $w(u(k))=\sum_{t \geq 1} \theta^k_t w(u^k_t)\leq w(p)$. So we   obtain:$$w(p) \geq w(u)\geq r(u) \geq v(p)-\varepsilon.$$
 \noindent This is true for all $\varepsilon$, so $w\geq v$.  \hfill $\Box$
 
 \begin{pro}  \label{pro5} $\Gamma$ has a general uniform value. 
  \end{pro}

\noindent{\bf Proof of proposition \ref{pro5}:}  First we can extend the notion of mixed play to $Z$. A mixed play at $u_0\in Z$, is a sequence $\sigma=(u_1,...,u_t,...)\in Z^\infty$ such that $u_{t+1} \in \hat{F}(u_{t})$ for each $t\geq 0$, and we denote by $\Sigma(u_0)$ the set of mixed play at $u_0$. Given $t, T$ in $\N$, $n \in \N^*$ and $u_0 \in Z$, we define for each mixed play $\sigma=(u_t)_{t \geq 1}\in \Sigma(u_0)$ the auxiliary payoff:
$${\gamma}_{t,n}(\sigma)= \frac{1}{n} \sum_{l=t+1}^{t+n} r(u_l), \; {\rm and } \;  \beta_{T,n}(\sigma)=\inf_{t \in \{0,...,T\}} {\gamma}_{t,n}(\sigma).$$

\noindent And we also define the auxiliary value function: for all $u$ in $Z$, $$h_{T,n}(u_0)=\sup_{\sigma\in \Sigma(u_0)} \beta_{T,n}(\sigma).$$
\noindent Clearly, $\beta_{T,n}(\sigma)\leq {\gamma}_{0,n}(\sigma)$ and $h_{T,n}(u_0)\leq v_n(u_0).$ We can write:
\begin{eqnarray*}
h_{T,n}(u_0) & = & \sup_{\sigma \in \Sigma(u_0)} \inf_{\theta\in \Delta(\{0,...,T\})}  \frac{1}{n}\sum_{t=0}^T \theta_t  \sum_{l=t+1}^{t+n} r(u_l)\\
& =&  \sup_{\sigma \in \Sigma(u_0)} \inf_{\theta\in \Delta(\{0,...,T\})} \sum_{l=1}^{T+n} \beta_l(\theta,n) r(u_l).
\end{eqnarray*}
 \noindent where for each $l$ in $1,...,T+n$, $$\beta_l(\theta,n)= \frac{1}{n}\sum_{t=Max\{0,l-n\}}^{\Min \{T,l-1\}}  \theta_t.$$
By construction, $\hat{F}$ is affine, so $\Sigma(u_0)$ is a convex subset of $Z^{\infty}$.  $ \Delta(\{0,...,T\}) $ is convex compact and the payoff  $\sum_{l=1}^{T+n} \beta_l(\theta,n) r(u_l)$ is affine both  in $\theta$ and   in $\sigma$.  We can apply a standard minmax theorem to get: 
$$h_{T,n}(u_0) = \inf_{\theta\in \Delta(\{0,...,T\})} \sup_{\sigma \in \Sigma(u_0)}\sum_{l=1}^{T+n} \beta_l(\theta,n) r(u_l).$$
We write $\theta_t=0$ for $t>T$ and  for each $l\geq 0$:
 $\beta_l(n,\theta)= \frac{1}{n} (\theta_0+...+ \theta_{l-1})$ if $l \leq n$,  $\beta_l(\theta,n)= \frac{1}{n} (\theta_{l-n}+...+ \theta_{l-1})$ if $n+1 \leq l \leq {n+T}$,   $\beta_l(n,\theta)=0$ if $l>n+T$. The evaluation  $\beta(\theta, n)$ is a particular probability on stages and $h_{T,n}(u_0)=  \inf_{\theta\in \Delta(\{0,...,T\})} v_{\beta(\theta,n)}(u_0).$
   It is easy to  bound the  impatience of $\beta(\theta,n)$: $$\sum_{l\geq 0} |\beta_{l+1}(\theta,n)-\beta_l(\theta,n)|= \sum_{l =0}^{n-1} \frac{\theta_l}{n} + \sum_{l \geq n} \frac{1}{n} | \theta_l - \theta_{l-n}| \leq \frac{3}{n}\longrightarrow_{n \to \infty} 0.$$
 
\noindent The impatience of $\beta(\theta,n)$ goes to zero as $n$ goes to infinity, uniformly in  $\theta$. So  we can use the previous proposition \ref{pro4}  to get: 
$$\forall \varepsilon>0, \exists n_0, \forall n \geq n_0, \forall \theta \in \Delta(\N), \forall u_0 \in Z, \; |v_{\beta(\theta,n)}(u_0)-v^*(u_0)| \leq \varepsilon.$$

This implies that  $h_{\infty,n}(u_0):=_{def}\inf_{\theta\in \Delta(\N)} v_{\beta(\theta,n)}(u_0)=\inf_{T\geq 0} h_{T,n}(u_0)$ converges to $v^*(u_0)$ when $n\to \infty$, and the convergence is uniform over $Z$. Consequently, if we fix  $\varepsilon>0$ there  exists $n_0$ such that for all $u_0$ in $Z$, for all  $T\geq 0$, there exists a play $\sigma^T=(u_t^T)_{t\geq 1}$ in $\Sigma(u_0)$ such that the average payoff is good on every interval of $n_0$ stages starting before  $T+1$:  for all $t=0,...,T$, $\; \gamma_{t,n_0}(\sigma^T)\geq v^*(u_0)-\varepsilon$.

We fix $u_0$ in $Z$ and consider, for each $T$, the play $\sigma^T=(u_t^T)_{t\geq 1}$ in $\Sigma(u)$ as above.  By a diagonal argument we can construct for each $t\geq 1$ a limit point $z_t$  in $\overline{Z}$  of the sequence $(u_t^T)_{T\geq 0}$ such that for each $t$ we have $(z_{t},z_{t+1})\in \cl (Graph (\hat{F}))$, with $z_0=u_0$. For each $m\geq 0$, we have $\frac{1}{n_0} \sum_{t=m+1}^{m+1+n_0} r(u_t^T)\geq v^*(u_0)-\varepsilon$ for $T$ large enough, so at the limit we get: $\frac{1}{n_0} \sum_{t=m+1}^{m+1+n_0} r(z_t)\geq v^*(u_0)-\varepsilon$.

$r$ being uniformly continuous, there exists $\alpha$ such that $|r(z)-r(z')|\leq \varepsilon$ as soon as $d(z,z')\leq \alpha$.   By lemma \ref{lem43}, one can find a $\sigma'=(z'_1,....,z'_t,...)$ at $\Sigma(z_0)$ such that for each $t$, $d(z_t,z'_t)\leq \alpha$. We obtain that for each $m\geq 0$,  $\frac{1}{n_0} \sum_{t=m+1}^{m+1+n_0} r(z'_t)\geq v^*(u)-2\varepsilon$. 

Consequently we have proved: $\forall \varepsilon>0$, there exists $n_0$ such that for all initial state $p$ in $X$, there exists a mixed play $\sigma'=(z'_t)_t$ at $p$ such that:  $\forall m\geq 0$,  $\frac{1}{n_0} \sum_{t=m+1}^{m+1+n_0} r(z'_t)\geq v^*(p)-2\varepsilon$.   Let $\theta\in \Delta(\N^*)$ be an evaluation,  it is now easy to conclude. First if $v^*(p)-2\epsilon<0$, then any play is $2\epsilon$-optimal. Otherwise, for each $j\geq 1$,   denote by $\overline{\theta_j}$ the maximum of $\theta$ on the block $B^j=\{(j-1) n_0+1,..., j n_0\}$. For all $t \in B^j$, we have:
\[\overline{\theta_j} \geq \theta_t  \geq \overline{\theta_j} - \sum_{t'\in \{(j-1) n_0+1,...j n_0-1\}} |\theta_{t'+1}-\theta_{t'} |.
\]
\noindent As a consequence, for all $j$ we have: \begin{align*}
\sum_{t=(j-1) n_0+1}^{j n_0}  \theta_t r(z'_t) & \geq   \overline{\theta_j}  \sum_{t=(j-1) n_0+1}^{j n_0}r(z'_t) \; - \; n_0 \sum_{t'\in \{(j-1) n_0+1,...,j n_0-1\}} |\theta_{t'+1}-\theta_{t'}| \\
\; & \geq    \sum_{t=(j-1) n_0+1}^{j n_0}\theta_t (v^*(p)-2\varepsilon) \; - \; n_0 \sum_{t'\in \{(j-1) n_0+1,...,j n_0-1\}} |\theta_{t'+1}-\theta_{t'}| \\
\end{align*}
and by summing over $j$, we get:  $\gamma_{\theta}(x_0,\sigma') \geq   v^*(p)-2\epsilon-n_0 I(\theta)\geq v^*(p)-3\epsilon$ as soon as $I(\theta)$ is small enough.  \hfill $\Box$

\subsection{Proof    of theorem \ref{thm4}} \label{sec35}
 
  Assume that $X$ is a compact subset of a simplex  $\Delta(K)$, and  
 let   $\Psi=(X,A,q,g)$ be a standard MDP such that:  $\forall x \in X, \forall y\in X, \forall a \in A,   \forall f \in D_1,  \forall \alpha \geq 0, \forall \beta \geq 0,$
$$ |\alpha f(q(x,a))-\beta f(q(y,a))|\leq \|\alpha x -\beta y\|_1 \; {\rm and}\; 
  |\alpha g(x,a)-\beta g(y,a)| \leq \|\alpha x -\beta y\|_1
.$$
We   write   $Z= \Delta_f(X)\times [0,1]$, and $\overline{Z}=\Delta(X) \times [0,1]$. We will use the metric $d_*=d_0=d_1=d_2=d_3$ on $\Delta(\Delta(K))$ introduced in section \ref{sub23} and its restriction to $\Delta(X)$, so that $\overline{Z}$ is a compact metric space.  For all $(u,y),(u',y')\in \Delta_f(X)\times[0,1]$, we put  $d((u,y),(u',y'))=\max(d_*(u,u'),|y-y'|)$ so that $(Z,d)$ is a precompact metric  space. Recall  we have defined the correspondence $\hat{F}$ from  $Z$ to itself such that for all $(u,y)$ in $Z$, $$\hat{F}(u,y)=\left\{(Q(u, \sigma),G(u, \sigma)) \; s.t. \ \sigma: X\rightarrow \Delta_f(A) \right\},$$

\noindent  with the notations $Q(u, \sigma)=\sum_{x \in X} u(x) q(x,\sigma(x))$ and $G(u , \sigma)= \sum_{x \in X} u(x) g(x,\sigma(x)).$
And we simply define the payoff function $r$ from $Z$ to $[0,1]$ by $r(u,y)=y$ for all $(u,y)$ in $Z$.  We start with a crucial lemma, which shows the importance of  the duality formula of theorem \ref{thm2}.  
 
   \begin{lem} \label{lem*} $\hat{F}$ is an affine and non expansive correspondence from $Z$ to itself. 
\end{lem}

\noindent{\bf Proof of lemma \ref{lem*}.} We first show that: $\forall u, u' \in \Delta_f(X)$, $\forall \alpha \in [0,1]$, $\forall y,y'\in [0,1]$,  $\hat{F}(\alpha u +(1-\alpha)u',\alpha y +(1-\alpha)y') = \alpha \hat{F}(u,y)+(1-\alpha) \hat{F}(u',y').$ First the transition does not depend on the second coordinate so we can forget it for the rest of the proof. The $\subset $ part is clear.  To see the reverse inclusion, consider  $\sigma: X\rightarrow \Delta_f(A)$, $\sigma':X\rightarrow \Delta_f(A)$ and   $v= \alpha \sum_{x \in X} u(x) q(x,\sigma(x)) + (1-\alpha) \sum_{x \in X} u'(x) q(x,\sigma'(x))$   in $\alpha \hat{F}(u)+(1-\alpha) \hat{F}(u')$.   Define $$\sigma^*(x)=\frac{\alpha u(x) \sigma(x) + (1 -\alpha) u'(x) \sigma'(x)}{\alpha  u(x)+(1-\alpha) u'(x)},$$ \noindent for each $x$ such that the denominator is positive. Then  $v=\sum_{x\in X} (\alpha u +(1-\alpha)u'(x))q(x,\sigma^*(x))  
$,  and $\hat{F}$ is affine. \\

We now prove that  $\hat{F}$ is non expansive. Let $z=(u,y)$ and  $z'=(u',y')$ be in $Z$. We have $d((u,y),(u',y'))\geq d_*(u,u')$ and denote by $U$ and $U'$ the respective supports of $u$ and $u'$. By the duality formula of theorem \ref{thm2}, there exists $\alpha=(\alpha(p,p'))_{(p,p')\in U \times U'}$  and $\beta=(\beta(p,p'))_{(p,p')\in U \times U'}$ with non-negative  coordinates satisfying: $\sum_{p' \in U'} \alpha(p,p')= u(p)$ for all $p\in U$, $\sum_{p \in U} \beta(p,p')=u'(p')$ for all $p'\in U'$, and 
     $$d_*(u,u')=\sum_{(p,p') \in U \times U'} \| p\;  \alpha(p,p')- p' \;  \beta(p,p')\|_1. $$
     
Consider now $v=Q(u, \sigma)=\sum_{p \in U} u(p) q(p,\sigma(p))$ for some  $\sigma:X\rightarrow \Delta_f(A)$. We define for all $p'$ in  $U'$:
\[
\sigma'(p')=\sum_{p\in U} \frac{\beta(p,p')}{u'(p')} \sigma(p),
\]
and $v'=Q(u', \sigma')= \sum_{p' \in U'} u'(p') q(p',\sigma'(p')).$
Then $v'\in \hat{F}(u',y')$, and for each test function $\varphi$ in $D_1$ we have:
\begin{align*}
|\varphi(v)-\varphi(v')| & =|\sum_{p,p'} \alpha(p,p') \varphi(q(p,\sigma(p))) - \beta(p,p') \varphi(q(p',\sigma(p))) | \\
 & =|\sum_{p,p',a} \alpha(p,p') \sigma(p)(a)\varphi(q(p,a)) - \beta(p,p') \sigma(p)(a)\varphi(q(p',a)) | \\
 & \leq \sum_{p,p'} \|\alpha(p,p') p -\beta(p,p') p' \|_1 = d_*(u,u'),
\end{align*}
\noindent and therefore $d_*(v,v')\leq d_*(u,u').$ In addition we have a similar result on the payoff, 
\begin{align*}
|G(u,\sigma)-G(u',\sigma')| &= |\sum_{p,p'} \alpha(p,p') g(p,\sigma(p)) - \beta(p,p') g(p',\sigma(p)) | \\
 & \leq \sum_{p,p'} \|\alpha(p,p') p -\beta(p,p') p' \|_1 \\
 & \leq d_*(u,u').
\end{align*}

\noindent Thus we have $d((Q(u, \sigma), R(u, \sigma)), (Q(u', \sigma'), R(u', \sigma'))) \leq d_*(u,u') \leq d(z,z')$. \hfill $\Box$

\vspace{0,5cm}

Recall that the set of invariant couples of the MDP $\Psi$ is:  $$RR=\{(u,y)\in \overline{Z}, ((u,y),(u,y))\in  cl(Graph(\hat{F}))\},$$ and the function $v^*:X \longrightarrow  \R$  is defined by: \begin{eqnarray*} v^*(x) & = \inf &\big\{ 
w(x),  w:\Delta(X)\rightarrow [0,1] \;  {\rm affine}\;  C^0 \; s.t. \\
 & \; &  (1)\;   \forall y\in X, w(y)\geq \sup_{a \in A} w(q(y,a))  \; {\rm and }\;      (2) \;  \forall (u,y) \in RR, w(u)\geq y \;  \big \}.\end{eqnarray*}

We now consider the deterministic Gambling House $\hat{\Gamma}=(Z, \hat{F}, r)$. $Z$ is precompact metric, $\hat{F}$ is affine non expansive and $r$ is obviously affine and uniformly continuous. Given an evaluation $\theta$,  the $\theta$-value  of $\hat{\Gamma}$ at $z_0=(u,y)$ is denoted by $\hat{v}_{\theta}(u,y)= \hat{v}_{\theta}(u)$ and does not depend on $y$.  The recursive formula of section \ref{sub31} yields:  \begin{eqnarray*}
\forall (u,y) \in Z,  \; \; \hat{v}_{\theta}(u)& =&\sup_{(u',y')\in \hat{F}(u)} \theta_1 y'+ (1- \theta_1) \hat{v}_{\theta^+} (u') \\
& = & \sup_{\sigma \in X \rightarrow \Delta_f(A)} \left( \theta_1 G(u,\sigma) + (1-\theta_1) \hat{v}_{\theta^+} (Q(u, \sigma))\right).
\end{eqnarray*}\noindent Because $\hat{F}$ and $r$ are affine, $\hat{v}_\theta$ is affine in $u$ and the supremum in the above expression can be taken over the function from $X$ to $A$. Because $\hat{F}$ is non expansive and $r$ is 1-Lipschitz,  each $\hat{v}_{\theta}$ is 1-Lipschitz. 

 We denote by $v_\theta$ the $\theta$-value of the MDP  $\Psi$ and linearly extend it to $\Delta_f(X)$. It turns out that the recursive formula satisfied by $v_\theta$ is similar to the above recursive formula for $\hat{v}_\theta$, so that $v_\theta(u)=\hat{v}_\theta(u,y)$ for all $u$ in $\Delta_f(X)$ and $y$ in $[0,1]$. As a consequence, the existence of the general limit value in both problems $\hat{\Gamma}$ and $\Psi$ is  equivalent. Moreover, a deterministic play in $\hat{\Gamma}$ induces a strategy in $\Psi$, so that the existence of a general uniform  value in $\hat{\Gamma}$  will imply  the existence of the  general uniform value in $\Psi$ (note that deterministic and mixed plays in $\hat{\Gamma}$ are equivalent since $\hat{F}$ has convex values). 
 
It is thus sufficient to show that $\hat{\Gamma}$ has a general uniform  value given by $v^*$, and we can   mimic the end of the proof of theorem \ref{thm3}.   Lemma \ref{lem43} applies word for word. Finally,  one can   proceed almost  exactly as in  propositions \ref{pro4} and \ref{pro5} to show that $\hat{\Gamma}$, hence $\Psi$,  has a general uniform value given by $v^*$.

 \section{Applications to partial observation and games} \label{sec4}

\subsection{POMDP with finitely many states} \label{sub42}

We now consider a more general model of MDP with actions where after each stage, the decision maker does not perfectly observe the  state.  A MDP with partial observation, or POMDP, $\Gamma=(K,A,S,q,g)$ is given by a finite set of states $K$, a non empty set of actions $A$ and a non empty set of signals $S$. The transition $q$ now goes from  $K \times A$ to $\Delta_f(S\times K)$ (by assumption the support of the signals at each state is finite) and the payoff function $g$ still goes from $K\times A$ to $[0,1]$. Given an initial probability $p_1$ on $K$, the POMDP $\Gamma(p_1)$ is played as following. An initial state $k_1$ in $K$ is selected according to $p_1$  and is not told to the decision maker. At every stage $t$ he selects an action $a_t\in A$. He has a (unobserved) payoff $g(k_t,a_t)$ and a pair $(s_t, k_{t+1})$ is drawn according to $q(k_t,a_t)$. The player learns $s_t$, and the play proceeds to stage $t+1$ with the new state $k_{t+1}$. A behavioral strategy is now  a sequence $(\sigma_t)_{t \geq 1}$ of applications with for each $t$, $\sigma_t:(A\times S)^{t-1}\rightarrow \Delta_f(A)$. As usual, an initial probability on $K$ and a behavior strategy $\sigma$ induce a probability distribution over $(K \times A \times S)^{\infty}$ and we can define the $\theta$-values and the notions of general limit and uniform values accordingly.

\begin{thm} \label{thm5} A POMDP with finitely many states  has a general uniform value, i.e. there exists $v^*:\Delta(K) \to \R$ with the following property:     for each $\varepsilon>0$ one can find $\alpha>0$ and for each initial probability $p$ a behavior  strategy $\sigma(p)$ such that for each evaluation $\theta$ with  $I(\theta) \leq \alpha$,
$$\forall p \in \Delta(K),  |v_{\theta}(p)-v^*(p)|\leq \varepsilon \; \; {and}\;\;   \gamma_{\theta}(\sigma(p))\geq v^*(p)-\varepsilon.  $$  \end{thm}

\noindent{\bf Proof:} We introduce $\Psi$ an auxiliary MDP on $X=\Delta(K)$ with the same set of actions $A$ and the following payoff and transition functions:
 \begin{quote} 
  $\bullet$   $r:X\times A \longrightarrow [0,1]$ such that $r(p,a)=\sum_{k\in K} p(k)g(k,a)$ for all $p$ in $X$ and $a\in A$,
  
  $\bullet$  $\hat{q}: X\times A \rightarrow \Delta_f(X)$ such that   
 $$\hat{q}(p,a)=\sum_{s \in S} \left( \sum_k p^k q(k,a)(s)\right) \delta_{\hat{q}(p,a,s)},$$
where $\hat{q}(p,a,s)\in \Delta(K)$ is the belief on the new state after playing $a$ at $p$ and observing the signal $s$:
$$\forall k'\in K, \hat{q}(p,a,s)(k')=\frac{q(p,a)(k',s)}{q(p,a)(s)}=\frac{\sum_k p^k q(k,a)(k',s)}{\sum_{k} p^k q(k,a)(s)}.$$

\end{quote}

The POMDP $\Gamma(p_1)$ and the standard MDP $\Psi(p_1)$  have the same value for all $\theta$-evaluations. And for each strategy $\sigma$ in $\Psi(p_1)$, the player can guarantee the same payoff in the original game  $\Gamma(p_1)$ by mimicking the strategy $\sigma$. So if we prove that $\Psi$ has a general uniform value it will imply that the POMDP $\Gamma$ has a general uniform value.


To conclude the proof, we will simply apply theorem  \ref{thm4} to the MDP $\Psi$.  We need to check the assumptions on the payoff and on the transition.

Consider any  $p$, $p'$ in 
$X$,   $a\in A$,   $\alpha\geq 0$ and $\beta \geq 0$.  We have:
\[
|\alpha r(p,a)-\beta r(p',a)|=|\sum_k (\alpha p(k)-\beta p'(k))g(k,a)| \leq \|\alpha p-\beta p' \| 
\]
Moreover for any   $f\in D_1$, we have:
\begin{align*}
| \alpha \hat{q}(p,a)(f)-\beta \hat{q}(p',a)(f) |& = |\sum_{s\in S} \left(\alpha q(p,a)(s) f(\hat{q}(p,a,s)) - \beta q(p',a)(s) f(\hat{q}(p',a,s)) \right)| \\
& \leq \sum_{s} \|\alpha q(p,a)(.,s) - \beta q(p',a)(.,s)\| \\
& \leq \sum_{s,k,k'} |\alpha p(k')q(k',a)(k,s)- \beta p'(k')q(k',a)(k,s) |\\
& \leq \sum_{s,k,k'} q(k',a)(k,s)|\alpha p(k')-\beta p'(k')|  =\|\alpha p-\beta p'\|.
\end{align*}
where the first  inequality comes from the definition of $D_1.$

By theorem \ref{thm4}, the MDP $\Psi$  has a general uniform value and we deduce that the  POMDP $\Gamma$ has a general uniform value. $\hfill \Box$

\begin{exa}\label{dark}
\rm Let $\Gamma=(K,A,S,q,g,p_1)$ be a POMDP where $K=\{k_1,k_2\}$, $A=\{a,b\}$, $S=\{s\}$ and $p_1=\delta_{k_1}$. The initial state is $k_1$ and since there is only one signal, the decision maker will obtain no additional information on the state. We say that he is in the dark. The payoff is given by $g(0,a)=g(0,b)=g(1,b)=0$ and $g(1,a)=1$, and the transition  by $q(1,a)=q(1,b)=\delta_{1,s}$, $q(0,a)=\delta_{0,s}$ and $q(0,b)=\frac{1}{2}\delta_{0,s}+\frac{1}{2}\delta_{1,s}.$  On one hand if the decision maker plays $a$ then the state stays the same and he receives a payoff of $1$ if and only if the state is $1$, on the other hand if he plays $b$ then he receives a payoff of $0$ but the probability to be in state $1$ increases. 

We define the function $r$ from $X=\Delta(K)$ to $[0,1]$ by $r((p,1-p),a)=1-p$ and $r((p,1-p),b)=0$ for all $p\in[0,1]$, and the function $\hat{q}$ from $X$ to $\Delta_f(X)$ by 
\[
\hat{q}((p,1-p),a)=\delta_{(p,1-p)} \text{ and }\hat{q}((p,1-p),b)=\delta_{(p/2,1-p/2)}.
\]
\noindent Then the standard MDP $\Psi=(\Delta(K),A,r,\hat{q})$ is the MDP associated in the previous proof to $\Gamma$. This MDP is deterministic since the decision maker is in the dark.

In this example, the existence of a general uniform value is immediate. If we fix $n\in \N$, the strategy $\sigma =b^n a^{\infty}$ which plays $n$ times $b$ and then $a$ for the rest of the game, guarantees a stage payoff of $(1-\frac{1}{2^n})$ from stage $n+1$ on, so the game has a general uniform value equal to $1$. Finally if we consider the discounted evaluations, one can show that the speed of convergence of $v_\lambda$ is slower than $\lambda:$
\[v_{\lambda}(p_1)=1-\frac{\ln(\lambda)}{\ln(2)}\lambda + O(\lambda).\]
All the spaces are finite but the partial observation implies that the speed of convergence is slower than $\lambda$ contrary to the perfect observation case where it is well known that the convergence is in $O(\lambda)$.
\end{exa}

\begin{rem}
\rm It is unknown if the  uniform value exists in pure strategies, i.e. if the behavior strategies $\sigma(p)$ of theorem \ref{thm5} can be chosen with values in $A$. This was already an open problem for the Cesàro-uniform value (see Rosenberg {\it  et al.}  2002 and Renault 2011 for different proofs requiring the use of behavioral strategies). In our proof, there are two related places where the use of lotteries on actions  is important. First in the proof of the convergence of the function $h_{T,n}$ (within the proof of theorem \ref{thm3}), we used Sion's theorem in order to inverse a supremum and an infimum so we need the convexity of the set of strategies. Secondly when we prove that the extended transition is $1$-Lipschitz (see lemma \ref{lem*}), the coupling between the two distributions $u$ and $u'$ introduces some randomization.
\end{rem}

\subsection{Zero-sum repeated games with an informed controller} \label{sub43}
 
We finally  consider  zero-sum repeated games with an informed controller. We start with  a general model $\Gamma=(K,I,J,C,D,q,g)$ of zero-sum repeated game, where we have 5 non empty finite sets:  a set of states $K$,  two  sets  of actions $I$ and $J$ and two   sets of signals $C$ and $D$, and we also have  a transition mapping   $q$ from  $K \times I \times J$ to $\Delta(K \times C\times D)$ and a payoff  function $g$ from $K\times I\times J$ to $[0,1]$.  Given an  initial probability $\pi$ on $\Delta(K\times C\times D)$,  the game $\Gamma(\pi)= \Gamma (K,I,J,C,D,q,g, \pi)$ is played as follows: at stage $1$, a triple $(k_1,c_1,d_1)$ is drawn according to $\pi$, player $1$ learns  $c_1$ and player $2$ learns  $d_1$. Then simultaneously player 1 chooses an action $i_1$ in $I$  and player 2 chooses an action  $j_1$ in $J$. Player $1$ gets a (unobserved) payoff $r(k_1,i_1,j_1)$ and player $2$ the opposite. Then a new triple $(k_2,c_2,d_2)$ is drawn accordingly to $q(k_1,i_1,j_1)$. Player $1$ observes $c_2$, player $2$ observes $d_2$ and the game proceeds to the next stage, etc...

A (behavioral) strategy for player $1$ is a sequence $\sigma=(\sigma_t)_{t\geq 1}$ where for each $t\geq 1$, $\sigma_t$ is a mapping from $(C\times I)^{t-1} \times C$ to $\Delta(I)$. Similarly a strategy for player $2$ is a sequence of mappings $\tau=(\tau_t)_{t\geq 1}$ where for each $t\geq 1$ ,$\tau_t$ is a mapping from $(D\times J)^{t-1} \times D$ to $\Delta(J)$. We denote respectively by $\Sigma$ and $\Tau$ the set of strategies of  player $1$ and    player $2$. An initial distribution $\pi$ and a couple of strategies $(\sigma,\tau)$ defines for each $t$ a probability on the possible histories up to stage $t$. And by Kolmogorov extension theorem,  it can be uniquely extended to  a   probability on the set of infinite histories $(K\times C\times D \times I\times J) ^{+\infty}$.

Given $\theta$ an evaluation function, we define the $\theta$-payoff of $(\sigma,\tau)$ in $\Gamma(\pi)$ as the expectation under $\P_{\pi,\sigma,\tau}$ of the payoff function,
\[
\gamma_{\theta}(\pi,\sigma,\tau)=\E_{\pi,\sigma,\tau}\left(\sum_t \theta_t \; r(k_t,i_t,j_t) \right).
\]

By Sion's theorem the game   $\gamma_{\theta}(\pi)$ has a value:  \[ v_{\theta}(\pi)=\max_{\sigma\in \Sigma} \min_{\tau \in \Tau} \gamma_{\theta}(\pi,\sigma,\tau)= \min_{\tau \in \Tau}\max_{\sigma \in \Sigma} \gamma_{\theta}(\pi,\sigma,\tau),	
\]
and we can define the general limit  value  as in the MDP framework. Note that we do not ask the convergence to be uniform for all $\pi$ in $\Delta(K \times C \times D)$, because we will later make some assumptions, in particular  on the initial distribution.

\begin{defi} \label{defi38} The repeated game $\Gamma(\pi)=(K,I,J,C,D,q,g, \pi)$ has a  general limit value $v^*(\pi)$ if $v_{\theta}(\pi)$ converges to $v^*(\pi)$ when $I(\theta)$ goes to zero, i.e.: $$\forall \varepsilon >0, \exists \alpha >0, \forall \theta, \;\; \left( \; I(\theta) \leq \alpha\implies \left( |v_{\theta}(\pi)-v^*(\pi)|\leq \varepsilon \right)\; \right).$$ \end{defi}

 And we can define a general uniform value by symmetrizing the definition for MDP. 

\begin{defi} The repeated game  $\Gamma(\pi)$ has a general uniform value if it has a general limit value $v^*$ and for each $\varepsilon>0$ one can find $\alpha>0$ and a couple of strategies $\sigma^*$ and $\tau^*$ such that for all evaluations $\theta$ with $I(\theta)\leq \alpha$: 
\[\forall \tau \in \Tau, \gamma_{\theta}(\pi,\sigma^*,\tau)\geq v^*(\pi)-\varepsilon  \; \; {\rm and} \; \; 
  \forall \sigma \in \Sigma, \gamma_{\theta}(\pi,\sigma,\tau^*)\leq v^*(\pi)+\varepsilon  .\] 
\end{defi}

\vspace{0,5cm}

 We   now focus on the case of a repeated game with an informed controller. We follow the definitions introduced in Renault (2012a). The first one concerns the information of the first player. We assume that he is always fully informed of the state and of the signal of the second player:

\begin{ass} \label{ass1}
There exist two mappings $\widetilde{k}:C\rightarrow K$  and $\widetilde{d}:C\rightarrow D$ such that, if $E$ denotes $\{(k,c,d)\in K\times C\times D,\ \widetilde{k}(c)=k,\ \widetilde{d}(c)=d \}$, we have:    $\forall (k,i,j)\in K\times I \times J$, $q(k,i,j)(E)=1,$ and $\pi(E)=1$. 
\end{ass}

Moreover we will assume that only player $1$  has  a meaningful influence on the transitions,  in the following sense. 

\begin{ass}\label{ass2}
The marginal of the transition on $K\times D$ is not influenced by player 2's action. For $k$ in $K$, $i$ in $I$ and $j$ in $J$, we denote by $\bar{q}(k,i)$ the marginal of $q(k,i,j)$ on $K \times D$. 
\end{ass}

The second player may influence the signal of the first player but he can not prevent him neither to learn the state nor to learn his own signal. Moreover he can not influence his own information, thus he has no influence on his beliefs about the state or about the beliefs of player $1$ about his beliefs. A repeated game satisfying assumptions \ref{ass1} and \ref{ass2} is called a repeated game with an informed controller.  It was proved in Renault (2012a)   that for such games the Cesàro-uniform value  exists and we  will generalize it here to the general uniform value.\\

\begin{exa}\label{aumann}
\rm We consider $\Gamma$ a zero-sum repeated game with incomplete information as studied by Aumann and Maschler (see reference from 1995). It is defined by a finite family $(G^k)_{k \in K}$ of  payoff matrices in $[0,1]^{I\times J}$ and $p\in \Delta(K)$ an initial probability. At the first stage, some state $k$ is selected according to $p$ and told to player 1 only. The second player knows the initial distribution $p$ but not the realization. Then the matrix game $G^k$ is repeated over and over.  At each stage the players observe past actions but not their payoff. Formally it is a zero-sum repeated game $\Gamma=(K,I,J,C,D,q,g)$ as defined previously, with $C= K\times I \times J$ and $D=I\times J$, and for all $(k,i,j)\in K\times I\times J$, $g(k,i,j)=G^k(i,j)$ and $q(k,i,j)=\delta_{k,(k,i,j),(i,j)}$.  For all $p\in\Delta(K)$, we denote by $\Gamma(p)$ the game where the initial probability $\pi \in \Delta(K\times C\times D)$ is given by $\pi=\sum_{k\in K} p(k)\delta_{k,(k,i_0,j_0),(i_0,j_0)}$ with $(i_0,j_0)\in I \times J$ fixed.
 
For each $n$, we denote by $v_n(p)$ the value of the $n$-stage game with initial probability $p$, where the payoff is the expected mean average of the  $n$ first stages. It is known that it satisfies the following dynamic programming formula:    
     $$v_{n}(p)=  \sup_{a \in \Delta(I)^K} \left( \frac{1}{n} r(p,a)+ \frac{n-1}{n} \sum_{k\in K,i\in I} p^k a^k(i) v_{n-1}(\hat{q}(p,a,i))\right).$$
where $p\in \Delta(K)$, $r(p,a)=\min_j (\sum_k  p^k G^k(a^k,j)) $ and $\hat{q}(p,a,i)$ is the conditional belief on $\Delta(K)$ given $p$, $a$, $i$:
$$\hat{q}(p,a,i)(k')=\frac{\sum_{k} p(k)a^k(i)q(k,i)(k')}{\sum_{k} p(k)a^k(i)}.$$
Starting from a belief $p$ about the state, if player $2$ observes action $i$ and knows that the distribution of actions of player $1$ is $a$, then he updates his beliefs to $\hat{q}(p,a,i)$. Aumann and Maschler have proved that the limit value exists and is characterized by 
$$v^*=\cav  f^*=\inf\{v: \Delta(K)\to [0,1], v\; {\rm concave}\; v \geq f^*\},
$$
where $f^*(p)=Val\left(\sum_k p^k G^k\right)$ for all $p\in \Delta(K)$. The function $f^*$ is the value of the game, called the non-revealing game, where player $1$ is forbidden to use his information. 
\end{exa}
  
\begin{thm} \label{thm7}
A zero-sum repeated game with an informed controller has a general uniform value.
\end{thm}

\vspace{0,5cm}

\noindent{\bf Proof of theorem \ref{thm7}:} Assume that $\Gamma(\pi)=(K,I,J,C,D,q,g, \pi)$ is a repeated game with an informed controller. 
The proof will consist of  5  steps. First we introduce  an auxiliary  standard Markov Decision Process $\Psi(\hat{\pi})$ on the state space $X=\Delta(K)$.  Then we show that for all evaluations $\theta$, the  repeated  game $\Gamma(\pi)$ and the MDP $\Psi(\hat{\pi})$ have the same $\theta$-value. In step 3  we check  that  the MDP satisfies the assumption of theorem \ref{thm4} so it has a general limit value and a general uniform value $v^*$. As a consequence the repeated game has a general limit value $v^*(\pi)$.  Then we prove that player $1$ can use an $\epsilon$-optimal strategy of the auxilliary MDP  in order to guarantee $v^*(\pi)-\epsilon$ in the original game. Finally we prove that Player $2$ can play by blocks in the repeated game in order to guarantee $v^*(\pi)+\epsilon$. And we obtain that $v^*(\pi)$ can be guaranteed by both players in the repeated game, so   it is the general uniform value of $\Gamma(\pi)$. \\

For every $P \in \Delta(K\times C\times D)$, we denote by $\overline{P}$ the marginal of $P$ on $K\times D$ and we put $\hat{P}=\psi_D(\overline{P})$ where  $\psi_D$ is the disintegration with respect to $D$ (recall proposition \ref{pro3}): for all $\mu \in \Delta(K\times D)$, $\psi_D(\mu)=\sum_{d\in D} \mu(d) \delta_{\mu(.|d)}.$\\

{\bf Step 1:} We put $X=\Delta(K)$ and  $A=\Delta(I)^K$ and  for every $p$ in $X$, $a$ in $A$ and $b$ in $\Delta(J)$, we define: 
\begin{align*}
r(p,a,b) &= \sum_{(k,i,j)\in K\times I\times J} p^k a^k(i)b(j)g(k,i,j) \in[0,1],\\
 {r}(p,a) &= \inf_{b\in \Delta(J)} r(p,a,b)=  \inf_{j\in J} r(p,a,j), \\
\overline{q}(p,a) &= \sum_{(k,i)\in K\times I} p^k a^k(i) \overline{q}(k,i) \in \Delta(K\times D), \\
\hat{q}(p,a) &= \psi_D(\overline{q}(p,a))=\sum_{d\in D} \overline{q}(p,a)(d)\delta_{\hat{q}(p,a,d)} \in \Delta_f(X).
 \end{align*}
Here  $\hat{q}(p,a,d)\in \Delta(K)$ is the belief of the second player on the new state after observing the signal $d$ and knowing that player 1 has played $a$ at $p$:
\[
\forall k'\in K, \hat{q}(p,a,d)(k')=\frac{\overline{q}(p,a)(k',d)}{\overline{q}(p,a)(d)}=\frac{\sum_k p^k q(k,a(k))(k',d)}{\sum_{k} p^k q(k,a(k))(d)}.
\]
We define the  auxiliary MDP $\Psi=(X,A,\hat{q}, r)$, and denote the $\theta$-value  in the MDP by $\hat{v}_{\theta}$. The MDP with initial state $\hat{\pi}$ has strong links with the repeated game $\Gamma(\pi)$.   \\

{\bf Step 2:} By proposition 4.23, part b)  in Renault 2012a), we have for all evaluations $\theta$ with finite support:
$$v_{\theta}(\pi)=\hat{v}_{\theta}(\hat{\pi}).$$
\noindent The proof relies on the same recursive  formula satisfied by $v$ and $\hat{v}$, and the equality can be easily extended to any evaluation $\theta$.
\[
\forall \theta \in \Delta(\N^*), \forall p \in X,    \; v_{\theta} (p)= \sup_{a \in A} \inf_{b\in B} \left(\;  \theta_1 r(p,a,b) + (1-\theta_1) v_{\theta^+} (\hat{q}(p,a))\; \right).
\]
where  $v_{\theta^+}$ is naturally  linearly extended to $\Delta_f(X)$. As a consequence if $\Psi(\hat{\pi})$ has a general limit value so does the repeated game $\Gamma(\pi)$.  \\

{\bf Step 3:}  Let us check that $\Psi$ satisfies the assumption of \ref{thm4}. Consider $p$, $p'$ in $X$,  $a$ in $A$,  and  $\alpha\geq 0$ and $\beta \geq 0$. We have:
\begin{eqnarray*}
|\alpha r(p,a)-\beta r(p',a)|& \leq & \sup_{b \in \Delta(J)} |\alpha r(p,a,b)-\beta r(p',a,b)| \\
\; & \leq & \sup_{b \in  \Delta(J)} |\sum_{k\in K} \alpha p^k g(k, a^k,b) - \beta p'^k g(k,a^k,b)|  \\
\; & \leq & \sup_{b \in  \Delta(J)} \sum_{k \in K} |\alpha p^k -\beta p'^k|= \|\alpha p-\beta p' \|_1 .
\end{eqnarray*}

\noindent Moreover, let $\varphi: \Delta(K) \longrightarrow \R$ be in $D_1$.
 \begin{eqnarray*}
|\alpha \varphi(\hat{q}(p,a)) -\beta \varphi(\hat{q}(p',a))|&= & \sum_{d \in D} \left( \alpha \bar{q}(p,a)(d) \varphi(\hat{q} (p,a,d))-  \beta \bar{q}(p',a)(d) \varphi(\hat{q} (p',a,d))\right) \\
\; & \leq &   \sum_{d \in D} \| \alpha \; \bar{q}(p,a)(d) \;  \hat{q} (p,a,d) -  \beta \; \bar{q}(p',a)(d) \;  \hat{q} (p',a,d)\|_1 \\
\; & \leq &  \sum_{d \in D} \| \alpha \; (\bar{q}(p,a)(k',d))_{k'}   -  \beta \; (\bar{q}(p',a)(k',d))_{k'} \|_1 \\
\; & \leq &  \sum_{d \in D} \sum_{k \in K} \| \alpha p^k \; (\bar{q}(k,a)(k',d))_{k'}   -  \beta p'^k  \; (\bar{q}(k,a)(k',d))_{k'} \|_1 \\
\; & \leq &  \sum_{d \in D} \sum_{k' \in K} \sum_{k \in K}\bar{q}(k,a)(k',d)  | \alpha p^k    -  \beta p'^k|  = \|\alpha p-\beta p' \|_1 .
\end{eqnarray*}

 So $\Psi=(X,A,\hat{q},r)$ has a general limit value and a general uniform value that we denote by $v^*$. As a consequence,  $\Gamma(\pi)$ has a general limit value $v^*(\pi)$. \\

{\bf Step 4:} Given $\varepsilon >0$, there exists  $\alpha>0$ and a  strategy $\sigma$  in the MDP $\Psi(\hat{\pi})$ such that the $\theta$-payoff in the MDP   is large: $\hat{\gamma}_{\theta}(\hat{\pi}, \sigma)\geq v^*(\pi)-\varepsilon$  whenever $I(\theta)\leq \alpha. $  Moreover if we look at the end of the proof of theorem \ref{thm4} we can choose  $\sigma$ to be induced by a deterministic play in the Gambling House $\hat{\Gamma}$ with state space $Z=\Delta_f(X)\times [0,1]$. As a consequence one can mimic  $\sigma$ to construct  a strategy $\sigma^*$ in the original repeated game $\Gamma(\pi)$ such that:   $\forall \tau \in \Tau, \gamma_{\theta}(\pi,  \sigma^*, \tau)\geq v^*(\pi)-\varepsilon$  whenever $I(\theta)\leq \alpha. $\\

{\bf Step 5:}
Finally we show that player $2$ can also guarantee the value $v^*$ in the repeated game $\Gamma$. Note that in the repeated  game  he can not compute the state variable in $\Delta(K)$ without knowing the strategy of player $1$. Nevertheless he has no influence on the transition function so playing independently by large blocks will be sufficient for him in order to guarantee $v^*(\pi)$. We use the following characterization of the value proved in Renault (2012a): 
\[
v^*(\pi)=\inf_n \sup_m v_{m,n}(\pi).
\]
where $v_{m,n}$ is the value of the game with payoff function  the Cesàro mean  of the stage payoffs between stages $m$ and $m+n.$ We proceed as in proposition 4.22 of Renault 2012a. 
Fix $n_0 \geq 1$, then we consider the strategy $\tau^*$ which for each $j\in \N$, plays optimally in the game with the evaluation the  Cesàro mean for the payoffs on the block of stages $B^j=\{n_0 (j-1)+1,...,n_0 j\}$. Since player $2$ does not influence the state it is well defined and this strategy guarantees $\sup_{t\geq 0} v_{t,n_0}(z)$ for the overall Cesàro mean.

Let $\theta\in \Delta(\N^*)$ and $\sigma$ be a strategy of player $1$. For each $j\geq 1$,   denote by $\underline{\theta_j}$ the minimum   of $\theta$ on the block $B^j=\{(j-1) n_0+1,..., j n_0\}$.  We have
\begin{align*}
\gamma_{\theta}(\pi,\sigma,\tau^*) & = \sum_{j=1}^{+\infty} \E_{\pi,\sigma,\tau^*} \left( \sum_{t=(j-1) n_0+1}^{j n_0} \theta_t  \; g(k_t,a_t,b_t) \right) \\
& \leq \sum_{j=1}^{+\infty} n_0 \; \underline{\theta}_j  \sup_{t\geq 0} v_{t,n_0}(\pi) + n_0 \sum_{t=1}^{+\infty} |\theta_{t+1}-\theta_{t}| \\
& \leq \sup_{t\geq 0} v_{t,n_0}(\pi) + n_0 I(\theta).
\end{align*}
Given $\epsilon$, there exists $n_0$ such that $\sup_{t\geq 0} v_{t,n_0}(\pi) \leq v^*(\pi)+\epsilon.$ Fix $\alpha=\frac{\epsilon}{n_0}$ and $\tau^*$ defined as before then for all $\theta$ such that $I(\theta) \leq \alpha$, we have
\[\sup_{\sigma \in \Sigma} \gamma_{\theta}(\pi,\sigma,\tau^*) \leq v^*(\pi)+ 2\epsilon, \]

\noindent and this concludes the proof of theorem \ref{thm7}.  \hfill $\Box$

%

%

\begin{exa}
\rm The computation of the value in two-player repeated game with incomplete information is a difficult problem as shown in the next example introduced in Renault (2006) and partially solved by H\"{o}rner {\it  et al.} (2010). The value exists by a theorem in Renault (2006) but the value has been computed only for some values of the parameters. The set of states is $K=\{k_1,k_2\}$, the set of actions of player $1$ is $I=\{T,B\}$, the set of actions of player $2$ is $J=\{L,R\}$, and the payoff is given by
\begin{center}
\begin{tabular}{cccccc}
 & $\begin{matrix} L  & R \end{matrix}$ & & & $\begin{matrix} L  & R \end{matrix}$ &\\
$\begin{matrix} T \\ B \end{matrix}$ &$\begin{pmatrix} 1 & 0 \\ 0 & 0  \end{pmatrix}$ & and & $\begin{matrix} T \\ B \end{matrix}$ &$\begin{pmatrix} 0 & 0 \\ 0 & 1  \end{pmatrix} $&.\\
 & $k_1$ & & & $k_2$ &
\end{tabular}
\end{center}
The evolution of the state does not depend on the actions: at each stage the state stays the same with probability $p$ and changes to the other state with probability $1-p$. At each stage, both players observe the past actions played but only player $1$ is informed of the current state (with previous notation $C=K\times I \times J$ and $D=I\times J$). For each $p\in [0,1]$, it defines a repeated game $\Gamma^p$. In the case $p=1$, the matrix is fixed for all the game thus it is a repeated game with incomplete information on one side {\it à la} Aumann Maschler (1995). For all other positive values of $p$, the process is ergodic so the limit value is constant, and it is sufficient to study the case $p\in [1/2,1)$ by symmetry of the problem. H\"{o}rner {\it  et al.} (2010) proved that if $p\in[1/2,2/3)$, then the value is $v_p=\frac{p}{4p-1}$. If $p\geq 2/3$, we do not know the value except for $p^*$, the solution of $9 p^3-12 p^2+ 6p-1=0$, where one has $v_p=\frac{p^*}{1-3p^*+6(p^*)^2}$.
\end{exa}

\vspace{1cm}
 
  \noindent \large \bf  Acknowledgements. \rm \small 
  The authors gratefully acknowledge  the support of the Agence Nationale de la Recherche, under grant ANR JEUDY, ANR-10-BLAN 0112, as well as the PEPS project Interactions  INS2I ``Propriétés des Jeux Stochastiques de Parité à Somme Nulle avec Signaux".

\section{References}

   

\noindent Araposthathis A., Borkar V., Fern\'{a}ndez-Gaucherand E., Ghosh M. and S. Marcus (1993): Discrete-time controlled Markov Processes with average cost criterion: a survey. SIAM Journal of Control and Optimization, 31, 282-344.\\

\noindent  Astr\"{o}m, K.J. (1965): Optimal Control of Markov Processes with Incomplete State
Information, J. Math. Anal. Appl., 10, 1965, 174-205. \\

\noindent  Ash, R.B. (1972): Real Analysis and Probability, Probability and Mathematical Statistics, Academic Press.\\%

\noindent  Aubin, J.P. (1977):  Applied Abstract Analysis. Wiley.\\%

\noindent  Aumann, R.J. and M. Maschler (1995): Repeated games with
 incomplete information.  With the collaboration of R. Stearns.
 {Cambridge, MA: MIT Press.}\\

\noindent Bellman, R. (1957): Dynamic programming, Princeton University Press.\\

\noindent Birkhoff, G.D. (1931): Proof of the ergodic theorem, Proceedings of the National Academy of Sciences of the United States of America, National Academy of Sciences. \\

\noindent Blackwell, D. (1962): Discrete dynamic programming. The Annals of  Mathematical Statistics, 33, 719-726. \\

\noindent Borkar, VS (2000): Average cost dynamic programming equations for controlled Markov chains with partial observation, SIAM Journal on Control and Optimization\\

		
\noindent Borkar , V.S. (2003):  Dynamic programming for ergodic control with partial observations, Stochastic Process. Appl.\\

\noindent Borkar, V.S. (2007) Dynamic programming for ergodic control of markov chains under partial observations: a correction, SIAM Journal on Control and Optimization.\\

\noindent Choquet, G. (1956): Existence et unicit{\'e} des repr{\'e}sentations int{\'e}grales au moyen des points extr{\'e}maux dans les c{\^o}nes convexes, S{\'e}minaire Bourbaki, 4, 33--47. \\


%
%
\noindent Denardo, E.V. and B.L. Fox (1968): Multichain Markov renewal programs,
SIAM Journal on Applied Mathematics

\noindent Dubins, L.E. and Savage, L.J. (1965): How to gamble if you must: Inequalities for stochastic processes, McGraw-Hill New York\\

\noindent Dudley, R. M: Real Analysis and Probability. Cambridge University Press, second edition 2002.\\


\noindent Filar, JA and Sznajder, R. (1992): Some comments on a theorem of Hardy and Littlewood, Journal of optimization theory and applications, 75, 201--208\\




\noindent Hern\'{a}ndez-Lerma, O. (1989): Adaptive Markov Control Processes. Springer-Verlag.\\



\noindent Hordijk, A. and Kallenberg, LCM (1979): Linear programming and Markov decision chains, Management Science, 352--362\\

\noindent H\"{o}rner J., D. Rosenberg, E. Solan and N. Vieille (2010): Markov Games with One-Sided Information. Operations Research, 58,  1107-1115. \\%

\noindent Kallenberg, LCM (1994): Survey of linear programming for standard and nonstandard Markovian control problems. Part I: Theory, Mathematical Methods of Operations Research, 1--42\\


\noindent Lehrer, E. and S. Sorin (1992): A uniform Tauberian Theorem in Dynamic Programming. Mathematics of Operations Research, 17, 303-307.\\


%




\noindent Maitra, A.P. and Sudderth, W.D. (1996): Discrete gambling and stochastic games,  Springer Verlag\\


\noindent Mertens, J.F. (1987): Repeated games, Proceedings of the International Congress of Mathematicians, 1528--1577. American Mathematical Society.\\

\noindent Mertens, J-F. and A. Neyman (1981):  Stochastic games. International Journal of Game Theory, 1, 39-64.\\

\noindent  Mertens, J-F., S. Sorin and S. Zamir (1994): {Repeated Games}.   Parts A, B and C.  CORE Discussion Papers  9420, 9421 and 9422.\\

\noindent Mertens, J-F. and S. Zamir (1985): Formulation of Bayesian analysis for games with incomplete information, International Journal of Game Theory, 14, 1--29.\\


\noindent Neyman, A.  (2008):  Existence of Optimal Strategies in Markov Games with Incomplete Information. International Journal of Game Theory, 581-596.\\


\noindent Renault, J. (2006): The value of Markov chain games with lack of information on one side. Mathematics of Operations Research, 31, 490-512.\\

\noindent Renault, J. (2011): Uniform value in Dynamic Programming.    Journal of the European Mathematical Society,  vol. 13, p.309-330. \\

\noindent Renault, J. (2012a): The value of   Repeated Games with an informed controller.    Mathematics of operations Research,  37: 154--179.
   \\

\noindent Renault, J. (2012b): General long-term values in Dynamic Programming, mimeo. \\

\noindent Rhenius, D. (1974): Incomplete information in Markovian decision models: The Annals of Statistics, 6, 1327--1334. \\

\noindent Rosenberg, D., Solan, E. and N. Vieille (2002):  Blackwell Optimality in Markov Decision Processes with Partial Observation. The Annals of Statisitics, 30, 1178-1193. \\

\noindent Rosenberg, D., Solan, E. and N. Vieille (2004): Stochastic games with a single controller and incomplete information. SIAM Journal on Control and Optimization, 43, 86-110. \\%

\noindent Runggaldier, W. J. and {\L}. Stettner (1991): On the construction of nearly optimal strategies for a general
              problem of control of partially observed diffusions, Stochastics Stochastics Rep.\\
		
\noindent Sawaragi Y. and  Yoshikawa T. (1970), Discrete-Time Markovian Decision Processes with Incomplete State Observations, Ann. Math. Stat., 41, 1970, 78-86. \\

\noindent Shapley, L.~S. (1953): Stochastic games. Proc. Nat. Acad. Sci., 39:1095--1100. \\




\noindent  Sorin, S.  (2002): A first course on Zero-Sum repeated games. MathÈmatiques et Applications, SMAI,  Springer.\\

\noindent  Villani, C. (2003), Topics in optimal transportation, Amer. Mat. Society , 58\\

\noindent  Von Neumann J. (1932), Proof of the quasi-ergodic hypothesis, Proc. Nat. Acad. Sci. U.S.A. 18,70 – 82.\\


 \end{document}